\documentclass[11pt]{article}
\usepackage{amssymb,amsmath,latexsym,amscd}

\newtheorem{example}[equation]{Example}
\newtheorem{theorem}[equation]{Theorem}
\newtheorem{corollary}[equation]{Corollary}
\newtheorem{lemma}[equation]{Lemma}
\newtheorem{proposition}[equation]{Proposition}
\newtheorem{definition}[equation]{Definition}
\newtheorem{remark}[equation]{Remark}
\numberwithin{equation}{section}
\newcommand{\K}{{k}}

\newcommand{\Q}{{\mathbb Q}}
\newcommand{\abG}{\left(\begin{array}{c}
a\\
b\end{array}\right)\ot \left(\begin{array}{c}
1\\
0\end{array}\right)}

\newcommand{\ghGbar}{\left(\begin{array}{c}
g\\
h\end{array}\right)\ot \left(\begin{array}{c}
0\\
1\end{array}\right)}

\def\wt{\widehat}
\newcommand{\pr}{\prime}
\newcommand{\whp}{\widehat{\p}}
\newcommand{\wh}{\widehat}
\newcommand{\Curr}{Curr\,}
\newcommand{\Ctd}{\hbox{\text{\rm Ctd}}}
\newcommand{\ib}{{\overline{\iota}}}
\newcommand{\ot}{\otimes}
\newcommand{\ol}{\overline}
\newcommand{\Z}{\mathbb{Z}}

\newcommand{\fg}{{\mathfrak g}}
\newcommand{\fE}{{\mathfrak E}}
\newcommand{\fR}{{\mathfrak R}}

\newcommand{\autfun}{{\bf Aut}}
\newcommand{\p}{\partial}
\newcommand{\pAS}{\partial_{\A\ot\widehat{\cS}}}
\newcommand{\ga}{\alpha}

\newcommand{\gs}{\sigma}
\newcommand{\eps}{\epsilon}
\newcommand{\proof}{{\bf Proof\ \ }}
\newcommand{\qed}{\hfill $\Box$}
\newcommand{\bm}{{\bf m}}

\newcommand{\bG}{{\bf G}}

\newcommand{\gd}{\partial}
\newcommand{\gl}{\lambda}
\newcommand{\cU}{{\mathcal U}}
\newcommand{\cC}{{\mathcal C}}
\newcommand{\F}{{\mathcal F}}

\newcommand{\A}{{\mathcal A}}
\newcommand{\B}{{\mathcal B}}
\newcommand{\End}{\hbox{End\,}}

\newcommand{\lb}[2]{\left\lbrack{#1}_\gl{#2}\right\rbrack}
\newcommand{\C}{{\mathbb C}}
\newcommand{\cL}{{\mathcal L}}
\newcommand{\R}{{\mathcal R}}
\newcommand{\cS}{{\mathcal S}}
\newcommand{\cT}{{\mathcal T}}
\newcommand{\even}{{\ol{0}}}
\newcommand{\odd}{{\ol{1}}}
\newcommand{\Aut}{{\rm Aut}}
\newcommand{\bAut}{{\bf Aut}}

\newcommand{\Hom}{{\rm Hom}}
\newcommand{\id}{{\rm id}}

\newcommand{\Span}{{\rm Span}}

\newcommand{\Gal}{{\mathcal Gal}}
\newcommand{\Loop}{{\mathcal L}}
\newcommand{\kalg}{{k-alg}}
\newcommand{\kalgder}{{k-\delta alg}}
\newcommand\limind{\mathop{\oalign{lim\cr
\hidewidth$\longrightarrow$\hidewidth\cr}}}
\newcommand\Spec{\text{\rm Spec}}
\newcommand\Out{\text{\rm Out}}
\newcommand\bOut{\text{\bf Out}}
\newcommand \bone{{\mathbf 1}}
\newcommand\Pic{\text{\rm Pic}}
\newcommand\os{\overset}
\newcommand\us{\underset}
\newcommand\et{\text{\rm \'et}}
\newcommand\ct{\text{\rm ct}}
\newcommand\fppf{\text{\rm fppf}}
\newcommand\q{\quad}
\newcommand{\SL}{{\mathbf{SL}}}
\newcommand{\bGL}{{\mathbf{GL}}}
\newcommand{\Zar}{{\rm Zar}}
\newcommand{\wtpsi}{\widetilde{\psi}}

\title{Differential Conformal Superalgebras and their Forms }

\author{ Victor Kac,$\hbox{}^{1 }$ Michael Lau,$\hbox{}^{2 *}$ and
Arturo Pianzola$\hbox{\,}^{3,4}$\thanks{Both authors gratefully
acknowledge the continuous support of the Natural Sciences and
Engineering Research Council of Canada.} \vspace{0.3cm}\\$\hbox{\
\,}^1${\small M.I.T., Department of Mathematics},\\ {\small
Cambridge, Massachusetts, USA 02139}\\ {\small Email:\
kac@math.mit.edu}\vspace{0.1cm}
\\$\hbox{\ \,}^2${\small University of Windsor,
Department of Mathematics and Statistics},\\ {\small Windsor,
Ontario, Canada N9B 3P4}\\ {\small Email:\ mlau@uwindsor.ca}\vspace{0.1cm}\\
$\hbox{\ \,}^3${\small University of Alberta, Department of
Mathematical and Statistical Sciences,}\\{\small Edmonton, Alberta,
Canada T6G 2G1}\\{\small Email:\ a.pianzola@math.ualberta.ca}\\
$\hbox{\ \,}^4${\small Instituto Argentino de Matem\'atica, Saavedra
15, (1083) Buenos Aires, Argentina}}
\date{}

\begin{document}
 \maketitle

\begin{small}
\noindent {\bf Abstract.} We introduce the formalism of differential
conformal superalgebras, which we show leads to the ``correct"
automorphism group functor and accompanying descent theory in the
conformal setting. As an application, we classify forms of $N=2$ and
$N=4$ conformal superalgebras by means of Galois cohomology.

\bigskip




\end{small}
\maketitle {\small \tableofcontents } \vskip14mm

 \setcounter{section}{-1} \vskip.25truein

\section{Introduction}
The families of superconformal algebras described in the work of
A.~Schwimmer and N.~Seiberg \cite{SchSe} bear a striking resemblance
to the loop realization of the affine Kac-Moody algebras \cite{K2}.
All of these algebras belong to a more general class known as {\em
$\Gamma$-twisted formal distribution superalgebras}, where $\Gamma$
is a subgroup of $\C$ containing $\Z$ \cite{FDA,KvdL}.

In a little more detail, a superconformal algebra\footnote{e.g., the Virasoro algebra or its superanalogues} (or
more generally, any twisted formal distribution algebra), is encoded by a
conformal superalgebra  $\A$ and an automorphism $\gs:\A\rightarrow\A$. Recall that $\A$ has a
$\C[\p]$-module structure and $n$-products $a_{(n)}b$, satisfying certain axioms \cite{kac}. Let
$\gs$ be a diagonalizable automorphism of $\A$ with eigenspace
decomposition
$$\A=\bigoplus_{\overline{m}\in\Gamma/\Z}\A_{\overline{m}},$$
where $\A_{\overline{m}}=\{a\in\A\ |\ \gs(a)=e^{2\pi i m}a\}$,
$\Gamma$ is an additive subgroup of $\C$ containing $\Z$, and
$\overline{m}\in\C/\Z$ is the coset $m+\Z\subset \C$.\footnote{  If
$\gs^M=1$, then this construction can be performed over an arbitrary
algebraically closed field $k$ of characteristic zero in the obvious
way, by letting $\Gamma$ be the group $\frac{1}{M}\Z$ and replacing
$e^{2\pi i/M}$ with a primitive $M$th root of $1$ in $k$. This is
the situation that will be considered in the present work.} Then the
associated $\Gamma$-twisted formal distribution superalgebra
$\hbox{Alg}\,(\A,\gs)$ is constructed as follows.

Let $\cL(\A,\gs)=\bigoplus_{m\in\Gamma}\left(\A_{\ol{m}}\ot_\C
t^m\right),$ and let
$$\hbox{Alg}\,(\A,\gs)=\cL(\A,\gs)\,\big/\,(\gd+\delta_t)\cL(\A,\gs),$$
where $\p$ denotes the map $\p\ot 1$ and $\delta_t$ is $1\ot\frac{d}{dt}$.  For each $a\in\A_{\ol{m}}$ and $m\in\Gamma$, let $a_m$ be the image
of the element $a\ot t^m\in \cL(\A,\gs)$ in $\hbox{Alg}\,(\A,\gs)$.
These elements span $\hbox{Alg}\,(\A,\gs)$, and there is a
well-defined product on this space, given by
\begin{equation}\label{0.1}
a_mb_n=\sum_{j\in\Z_+}{{m}\choose{j}}\left(a_{(j)}b\right)_{m+n-j},
\end{equation}
for all $a\in\A_{\overline{m}}$ and $b\in\A_{\overline{n}}$.

The name {\em twisted formal distribution algebra} comes from the fact that the
superalgebra $\hbox{Alg}(\A,\gs)$ is
spanned by the coefficients of the family of twisted pairwise local
formal distributions
$$\mathcal{F}=
\bigcup_{\overline{m}\in\Gamma/\Z}\left\{a(z)=
\sum_{k\in\overline{m}}a_kz^{-k-1}\ |\
a\in\A_{\overline{m}}\right\}.$$ For $\gs =1$ and $\Gamma = \Z$, we
recover the maximal non-twisted formal distribution superalgebra
associated with the conformal superalgebra $\A$.  See \cite{kac},\cite{FDA} for
details.

For example, let $A$ be an ordinary superalgebra over $\C$.  The {current conformal superalgebra}
$\A=\C[\p]\ot_\C A$ is defined by letting $a_{(n)}b=\delta_{n,0}ab$ for
$a,b\in A$ and extending these $n$-products to $\A$ using the conformal superalgebra
axioms.  The
associated {\em loop algebra} $A\ot_\C\C[t,t^{-1}]$ is then encoded by the current
superconformal algebra $\A$.  Taking $\gs$ to be an automorphism of $\A$ extended from a
finite order (or, more generally, semisimple) automorphism of $A$,
we recover the construction of a $\gs$-twisted loop algebra
associated to the pair $(\A , \gs)$. When $A$ is a Lie algebra, this
is precisely the construction of $\gs$-twisted loop algebras
described in \cite{K2}.

Under the correspondence described above, the superconformal
algebras on Schwimmer and Seiberg's lists are the $\Gamma$-twisted
formal distribution algebras associated with the $N=2$ and $N=4$ Lie
conformal superalgebras \cite{kac,KvdL}.  Prior to Schwimmer and
Seiberg's work, it was generally assumed that the $N=2$ family of
superconformal algebras consisted of infinitely many distinct
isomorphism classes.  However, it was later recognized that this
family contains (at most) two distinct isomorphism classes.  A
similar construction with $N=4$ superconformal algebras was believed
to yield an infinite family of distinct isomorphism classes
\cite{SchSe}.


The connection of the construction of the superalgebra
$\hbox{Alg}\,(\A,\gs)$ to the theory of differential conformal
superalgebras is as follows.  The $\C[\p]$-module $\cL(\A,\gs)$
carries the structure of a differential conformal superalgebra with
derivation $\delta=\delta_t$ and $n$-products given by
\begin{equation}
(a\ot t^k)_{(n)}(b\ot t^\ell)=\sum_{j\in\Z_+}{k\choose j} \big(a_{(n+j)}b\big)\ot t^{k+\ell-j}.
\end{equation}
Then $(\p+\delta)\cL(\A,\gs)$ is an ideal of $\cL(\A,\gs)$ with respect
to the $0$-product, which induces the product given by (\ref{0.1})
on $\hbox{Alg}\,(\A,\gs)$.  Moreover, the differential conformal
superalgebra $\cL(\A,\gs)$  is a twisted form of the affinization
$L(\A)=\cL(\A,id)$ of $\A$.

Thus, there are two steps to the classification of twisted superconformal algebras $\hbox{Alg}\,(\A,\gs)$ up to isomorphism.
First, we classify the twisted forms of the differential conformal superalgebra $L(\A)$.  In light of the above discussion, this
gives a complete (but possibly redundant) list of superconformal algebras, obtained by factoring by the image
of $\p+\delta$ and retaining only the $0$-product. Second, we should figure out which of these resulting superconformal algebras are
non-isomorphic.

The second step of the classification is rather straightforward. For
example, the twisted $N=4$ superconformal algebras are distinguished
by the eigenvalues of the Virasoro operator $L_0$ on the odd part.
The remainder of the paper will consider the first step, namely the
classification of the $\cL(\A, \gs)$ up to isomorphism.
\bigskip

Recently, the classification (up to isomorphism) of affine Kac-Moody
algebras has been given in terms of torsors and non-abelian \'etale
cohomology \cite{pianzola}.  The present paper develops conformal
analogues of these techniques, and lays the foundation for a
classification of forms of conformal superalgebras by cohomological
methods. These general results are then applied to classify the
twisted $N=2$ and $N=4$ conformal superalgebras up to isomorphism.

 To illustrate our methods, let us look at the case of the twisted loop
algebras  as
 they appear in the theory of affine Kac-Moody Lie algebras.
Any such  $\cL$  is
 naturally a Lie algebra over $R:=\C [t^{\pm 1}]$
and
\begin{equation}\label{split}
\cL \otimes _R S\simeq \fg \otimes _\C S \simeq (\fg\otimes _\C
R)\otimes _R S
\end{equation}
for some (unique) finite-dimensional simple Lie algebra $\fg,$ and
some (finite, in this case) \'etale extension $S/R.$  In particular,
$\cL$ is an $S/R$-form of the $R$-algebra $\fg\otimes_\C R$, with
respect to the \'etale topology of $\Spec(R)$.  Thus $\cL$
corresponds to a torsor over $\Spec(R)$ under $\autfun(\fg)$ whose
isomorphism class is an element of the pointed set $H^1_{\text{\rm
\'et}}\big(R,\autfun(\fg)\big).$

Similar considerations apply to forms of the $R$-algebra $A\otimes
_k R$ for any finite-dimensional
 algebra $A$ over an algebraically
closed field $k$ of characteristic $0.$ The crucial point in the
classification of forms of $A\otimes_k R$ by cohomological methods
is that in the exact sequence of pointed sets
\begin{equation}\label{exact}
H_{\et}^1 \big(R,\autfun^0(A)\big) \to H_{\et}^1\big(R,\autfun
(A)\big) \os\psi  \longrightarrow H_{\et}^1\big(R,\bOut(A)\big),
\end{equation}
 where $\bOut(A)$ is
the (finite constant) group of connected components of $A,$ the map
$\psi$ is injective \cite{pianzola}.

Grothendieck's theory of the algebraic fundamental group allows us
to identify $H_{\et}^1\big(R,\bOut(A)\big)$ with the set of
conjugacy classes of the corresponding finite (abstract) group
$\Out(A).$ The injectivity of the map $$H_{\et}^1\big(R,\autfun
(A)\big) \os\psi \longrightarrow H_{\et}^1\big(R,\bOut(A)\big)$$
means that to any form $\cL$ of $A \ot_k R,$  we can attach a
conjugacy class of the finite group $\Out(A)$ that characterizes
$\cL$ up to $R$-isomorphism. In particular, if $\autfun(A)$ is
connected, then all forms (and consequently, all twisted loop
algebras) of $A$ are {\it trivial}--that is, isomorphic to
$A\otimes_k R$ as $R$-algebras.

With the previous discussion as motivation, we now consider the
$N=2,4$ Lie conformal superalgebras $\A$ described in \cite{kac}.
 The automorphism groups of these objects are as follows:

\vskip.3cm
\centerline{Table 1}
\vskip.3cm \centerline{
\begin{tabular}{c|c|}
$N$ & $\Aut(\A)$\\
\hline
2 &$ \C^\times \rtimes \Z/2\Z $\\
\hline
4 &$\big(\SL_2(\C) \times \SL_2(\C)\big)/\pm I$\\
\hline
\end{tabular} }

\vskip.3cm It was originally believed that the standard $N=2$
algebra lead to an infinite family of non-isomorphic superconformal
algebras (arising as $\Gamma$-twisted formal distribution algebras
of the different $\cL(\A, \gs),$ as we explained above).
 This is somewhat surprising, for since $\C^\times \rtimes
\Z/2\Z$ has two connected components, one would expect (by analogy
with the finite-dimensional case) that there would be only two
non-isomorphic twisted loop algebras attached to $\A.$ Indeed,
Schwimmer and Seiberg later observed that all of the superconformal
algebras in one of these infinite families are isomorphic
\cite{SchSe}, and that (at most) two such isomorphism
classes existed.

On the other hand, since the automorphism group of the $N=4$
conformal superalgebra is connected, one would expect all twisted
loop algebras in this case to be trivial and, a fortiori, that all
resulting superconformal algebras would be isomorphic. Yet Schwimmer
and Seiberg aver in this case the existence of an infinite family of
non-isomorphic superconformal algebras!

The explanation of how, in the case of conformal superalgebras, a
connected automorphism group allows for an infinite number of
non-isomorphic loop algebras is perhaps the most striking
consequence of our work. Briefly speaking, the crucial point is as
follows.  A twisted loop algebra $\cL$ of a $k$-algebra $A$ is
always split by an extension $S_m:=k[t^{\pm1/m}]$ of $R:=k[t^{\pm
1}]$, for some positive integer $m$.  The extension $S_m/R$ is
Galois, and its Galois group can be identified with $\Z/m\Z$ by
fixing a primitive $m$th root of 1 in $k.$  The cohomology class
corresponding to $\cL$ can be computed using the usual Galois
cohomology $ H^1\big(\Gal(S_m/R),\autfun(A)(S_m)\big)$, where
$\autfun(A)(S_m) $ is the automorphism group of the $S_m$-algebra
$A\otimes_k S_m.$ One can deal with all loop algebras at once by
considering the direct limit $\widehat{S}$ of the $S_m $, which plays the
role of the ``separable closure" of $R.$ In fact, $\widehat S$ is the
simply connected cover of $R$ (in the algebraic sense), and the
algebraic fundamental group $\pi_1(R)$ of $R$ at its generic point
can thus be identified with $\widehat{\Z}$, namely the inverse limit
of the groups $\Gal(S_m/R) = \Z/m\Z .$

Finding the ``correct'' definitions of conformal superalgebras over
rings and of their automorphisms leads to an explanation of how
Schwimmer and Seiberg's infinite series in the $N = 4$ case is
possible. In our framework, rings are replaced by rings equipped
with a $k$-linear derivation (differential $k$-rings). The resulting
concept of differential conformal superalgebra is central to our
work, and one is forced to rewrite all the faithfully flat descent
formalism in this setting. Under some natural finiteness
conditions, we recover the situation that one encounters in the
classical theory, namely that the isomorphism classes of twisted
loop algebras of $\A$ are parametrized by $
H^1\big(\widehat{\Z},\autfun(\A)(\widehat{S})\big)$ with
$\widehat{\Z} = \Gal(\widehat{S}/R)$ acting continuously as
automorphism of $\A \otimes_k \widehat{S}.$

In the $N=2$ case, the automorphism group $\autfun (\A)(\widehat S)
= \widehat{S}^\times \rtimes \Z/2\Z$,  and  the cohomology set
$H^1\big(\widehat{\Z},\autfun(\A)(\widehat{S})\big) \simeq \Z/2\Z$, as
expected. By contrast, in the $N=4$ case, $\autfun(\A)(\widehat S)$
is {\em not} $ \big(\SL_2(\widehat S)\times \SL_2(\widehat
S)\big)/\pm I$ as we would expect from Table 1 above. In fact,
$$\autfun(\A)(\widehat S) = \big(\SL_2(\widehat S)\times \SL_2(\C)\big)/\pm I.$$
The relevant $H^1$ vanishes for  $\SL_2(\widehat S)$, but it is the
somehow surprising appearance of the ``constant" infinite group
$\SL_2(\C)$, through the (trivial) action of $\pi_1(R) = \widehat{\Z},$ that is ultimately responsible for an infinite family of non-isomorphic twisted conformal superalgebras that are
parametrized by the conjugacy classes of elements of finite order of
${\rm \bf PGL}_2(\C).$

In this paper, we use differential conformal (super)algebras for the
study of forms of conformal (super)algebras.  However, the theory
of differential conformal (super)algebras reaches far beyond.  For
example, it is an adequate tool in the study of differential (super)algebras; see Remark 2.7d in \cite{kac}. The ordinary conformal
(super)algebras do not quite serve this purpose since they allow
only translationally invariant differential (super)algebras.
Another area of applicability of differential conformal (super)algebras is
the theory of (not necessarily conservative) evolution PDEs.

\bigskip

\noindent {\bf Notation:} Throughout this paper, $k$ will be a field
of characteristic zero.  We will denote by $\kalg$  the category of
unital commutative associative $k$-algebras. If $k$ is algebraically
closed, we also fix a primitive $m$th root of unity $\xi_m\in k$ for
each $m>0$.  We assume that these roots of unity are {\em
compatible}.  That is, $\xi_{\ell m}^\ell=\xi_m$ for all positive
integers $\ell$ and $m$.

The integers, nonnegative integers, and rationals will be denoted
$\Z$, $\Z_+$, and ${\mathbb{Q}}$, respectively.  For pairs $a,b$ of
elements in a superalgebra, we let $p(a,b)=(-1)^{p(a)p(b)}$ where
$p(a)$ (resp., $p(b)$) is the parity of $a$ (resp., $b$).

Finally, for any linear transformation $T$ of a given $k$-space $V$,
and for any nonnegative integer $n$, we follow the usual convention
for divided powers and define $T^{(n)}:=\frac{1}{n!}T^n.$

\bigskip

\noindent {\bf Acknowledgments.}  Most of this work was completed
while M.L. was an NSERC postdoc at University of Ottawa and a
visiting fellow at University of Alberta.  He thanks  both
universities for their hospitality.  M.L. also thanks J.~Fuchs,
T.~Quella, and Z.~\v Skoda for helpful conversations.

\section{Conformal superalgebras}\label{egy}

This section contains basic definitions and results about conformal
superalgebras over differential rings.  We recall that $k$ denotes a
field of characteristic $0,$ and $\kalg$ is the category of unital
commutative associative $k$-algebras.

\subsection{Differential rings}

  For capturing the right concept of
conformal superalgebras over rings, each object in the category of
base rings should come equipped with a derivation. This leads us to
consider the category $\kalgder$ whose objects are pairs $\R =
(R,\delta_R)$ consisting of an object $R$ of $\kalg$ together with a
$k$-linear derivation $\delta_R$ of $R$ (a differential $k$--ring).
A morphism from $\R = (R,\delta_R)$ to $\cS = (S,\delta_S)$  is a
$k$-algebra homomorphism $\tau:\ R\rightarrow S$ that commutes with
the respective derivations.  That is, the diagram
\begin{equation}\label{ringextdiag}
\begin{CD}
R@> \tau>>S \\
@V\delta_RVV @VV\delta_SV \\
R@> \tau>>S.
\end{CD}
\end{equation}

\noindent commutes.

For a fixed $\R = (R,\delta_R)$ as above, the collection of all $\cS
= (S, \delta_S)$ in $\kalgder$ satisfying (\ref{ringextdiag}) leads
to a subcategory of $\kalgder$, which we denote by $\R-ext$. The
objects of this subcategory are called {\it extensions} of $\R$.
Each extension $(S,\delta_S)$ admits an obvious $R$-algebra
structure: $s\cdot r=r\cdot s:=\tau(r)s$ for all $r\in R$ and $s\in
S$.  The morphisms in $\R-ext$ are the $R$-algebra homorphisms
commuting with derivations. That is, for any $S_1,S_2\in\R-ext$,
$\Hom_{\R-ext}(S_1,S_2)$ is the set of $R$-algebra homomorphisms in
$\Hom_{\kalgder}(S_1,S_2)$.

 Let $\cS_i=\{(S_i,\delta  _i)\ |\ 1\le i\le n\}$ be a family
of extensions of $\R = (R,\delta  _R).$  Then
$$\delta:=\os
n{\us{i=1}\sum} \,\id \otimes \cdots \otimes \delta  _i\otimes
\cdots \otimes \,\id$$ is a $k$-linear derivation of $S_1\otimes_R\,
S_2\otimes_R \cdots \otimes_R\, S_n. $ The resulting extension
$(S_1\otimes_R\, S_2\otimes_R \cdots \otimes_R\, S_n, \delta)$ of
$(R, \delta_R)$ is called the {\it tensor product} of the $\cS_i$
and is denoted by $\cS _1\otimes_{\R} \cdots \otimes_{\R} \,\cS_n.$

Similarly, we define the {\it direct product} $\cS_1\times \cdots
\times \cS_n$ by considering the $k$-derivation $\delta _1 \times
\cdots \times \delta  _n$ of $S_1\times \cdots \times S_n.$

\begin{example}\label{Laurentextension} {\rm Consider the Laurent
polynomial ring $R = k[t,t^{-1}].$  For each positive integer $m$,
we set $S_m = k[t^{1/m}, t^{-1/m}]$ and
 $\widehat S = \limind S_m$.\footnote{In \cite{GiPi1}, \cite{GiPi2}, and \cite{GiPi3},
 where the multivariable case is considered, the rings
 $R$, $S_m$ and $\widehat S$ were denoted by $R_1$, $R_{1,m}$ and
 $R_{1,\infty}$ respectively.} We can think of $\widehat S$
 as the ring $k[t^q\,|\,q\in {\mathbb{Q}}]$ spanned by all rational powers
 of the variable $t$.
The $k$-linear derivation $\delta_t = \frac{d}{dt}$ of $R$ is also a
derivation of $S_m$ and $\widehat S$. Thus $\R = (R, \delta_t)$,
$\cS_m = (S_m , \delta_t)$, and $\widehat{\cS} = (\widehat S,
\delta_t)$ are objects in $\kalgder.$ Clearly $\cS_m$ is
 an extension of $\R$, and $\widehat{\cS}$ is an extension of $\cS_m$
 (hence also of $\R$). These rings with derivations will play a crucial
 role in
our work.}

\end{example}

\subsection{Differential conformal superalgebras}

Throughout this section,  $\R = (R, \delta_R)$ will denote an object
of $\kalgder$.

\begin{definition}\label{scdefinition}{\em {\em An $\R$-conformal
superalgebra} is a triple \\ $\big(\A, \p_{\A}, (\,-\, _{(n)} \,-\,
)_{n \in \Z_+}\big)$ consisting of
\medskip

(i) a $\Z/2\Z$-graded $R$-module $\A=\A_\even\oplus\A_\odd$,

(ii) an element $\p_{\A} \in \End_k(\A)$ stabilizing the even and
odd parts of $\A,$

(iii) a $k$-bilinear product $(a,b)\ \mapsto\ a_{(n)}b$ for each
$n\in\Z_+$,
\medskip

\noindent satisfying the following axioms for all $r\in R$,
$a,b,c\in\A$ and $m,n\in\Z_+$:

\begin{enumerate}
\item[{\rm (CS0)}] $a_{(n)}b=0$ for $n\gg 0$
\item[{\rm (CS1)}] $\p_\A (a)_{(n)}b=-na_{(n-1)}b$ and
$a_{(n)}\big(\p_\A (b)\big)=\p_\A\big(a_{(n)}b\big)+na_{(n-1)}b$
\item[{\rm (CS2)}] $\p_{\A}(ra) = r\p_{\A}(a)+\delta_R(r)a$
\item[{\rm (CS3)}] $a_{(n)}(rb)=r(a_{(n)}b)$ and
$(ra)_{(n)}b=\sum_{j\in\Z_+}\delta_R^{(j)}(r)\big(a_{(n+j)}b\big)$.
\end{enumerate}
\noindent If $n=0$, (CS1) should be interpreted as $\p_\A
(a)_{(0)}b=0$.
 Note that $\p_{\A}$ is a derivation of all $n$--products, called the {\it derivation of} $\A.$ The binary
 operation $(a,b)\ \mapsto\ a_{(n)}b$ is called the $n${\it -product of }
 $\A$.

If the $\R$-conformal superalgebra $\A$ also satisfies the following
two axioms, $\A$ is said to be an {\em $\R$-Lie conformal
superalgebra}:
\begin{enumerate}
\item[{\rm (CS4)}]$a_{(n)}b=-p(a,b)\sum_{j=0}^\infty(-1)^{j+n}
\p_{\A}^{(j)}(b_{(n+j)}a)$
\item[{\rm (CS5)}]$a_{(m)}(b_{(n)}c)=\sum_{j=0}^m\binom{m}{j}
(a_{(j)}b)_{(m+n-j)}c+p(a,b)b_{(n)}(a_{(m)}c).$
\end{enumerate}}
\end{definition}
\medskip

\begin{remark}\label{homothety} {\rm For a given $r \in R$ we will denote the corresponding
homothety $ a \mapsto ra$ by $r_{\A}.$ Axiom (CS2) can then be
rewritten as follows:
\begin{enumerate}
\item[{\rm (CS2)}]  $\p_{\A}\circ r_{\A} = r_{\A} \circ
\p_{\A}+ \delta_R(r)_{\A}$
\end{enumerate}
}
 \end{remark}

\begin{remark}\label{scoverfields} {\rm If $R = k$, then $\delta_R$
is necessarily the zero derivation.  Axioms (CS2) and (CS3) are then
superfluous, as they simply say that $\p_{\A}$  and that the
$(n)$-products are $k$-linear.  The above definition thus
specializes to the usual definition of conformal superalgebra over
fields (cf.~\cite{kac} in the case of complex numbers). Henceforth
when referring to a conformal superalgebra over $k$, it will always
be understood that $k$ comes equipped with the trivial derivation.}
\end{remark}

\begin{example}\label{loopexample} (Affinization of a conformal superalgebra) {\em
Let $\A$ be a conformal superalgebra over $k$.  As in Kac
\cite{kac},\footnote{The definition given in \cite{kac} is an
adaptation of affinization of vertex algebras, as defined by
Borcherds \cite{borcherds}.} we define the affinization
$L(\A)$ of $\A$ as follows. The underlying space of
$L(\A)$ is $\A\ot_{k} k[t,t^{-1}]$, with the $\Z/2\Z$-grading
given by assigning even parity to the indeterminate $t$ .  The
derivation of $L(\A)$ is
$$\gd_{L(\A)}=\p_\A\ot 1 + 1\ot \delta_t$$
where $\delta_t = \frac{d}{dt}$, and the $n$-product is given by
\begin{equation}\label{affprod}
(a\ot f)_{(n)}(b\ot
g)=\sum_{j\in\Z_+}(a_{(n+j)}b)\ot\delta_t^{(j)}(f)g
\end{equation}
for all $n\in\Z_+$, $a,b\in\A$, and $f,g\in k[t,t^{-1}]$. It is
immediate to verify that $\A$ is a $k$-conformal superalgebra. In
fact, $\A$ is in the obvious way a $(k[t,t^{-1}],
\delta_t)$--conformal superalgebra.
\medskip

If $R =k[t,t^{-1}]$  and  $\R = (R,\delta_t)$ are as in Example
\ref{Laurentextension}, then the affinization $L(\A)=\A\ot_k
R$ also admits an $\R$-conformal structure via the natural action of
$R$ given by $ r'(a\ot r):= a\ot r'r \label{i}$  for all $a\in\A$
and $r,r'\in R$.  The only point that needs verification is Axiom
(CS2),
 and this is straightforward to check.

Thus  $\A \otimes_k R$ is both a $k$- and an $\R$-conformal
superalgebra. We will need both of these structures in what follows.
From a physics point of view, it is the $k$-conformal structure that
matters; from a cohomological point of view, the $\R$-conformal
structure is crucial.

We will also refer to the affinization $L(\A)=\A\ot_k R$ as
the {\it (untwisted) loop algebra of} $\A$. It will always be made
explicit whether $L(\A)$ is being viewed as a $k$- or as an
$\R$-conformal superalgebra.}
\end{example}

Let  $\A$ and $\B$ be $\R = (R,\delta_R)$-conformal superalgebras. A
map $\phi:\ \A\rightarrow\B$ is called a {\em homomorphism} of
$\R$-conformal superalgebras if it is an $R$-module homomorphism
that is homogeneous of degree $\ol{0}$, respects the $n$-products,
and commutes with the action of the respective derivations. That is,
it satisfies the following three properties:
\begin{itemize}
\item[{\rm (H0)}] $\phi$ is $R$-linear and
$\phi(\A_{\ib})\subseteq\B_\ib$ for $\ib=\ol{0},\ol{1}$
\item[{\rm (H1)}] $\phi(a_{(n)}b)=\phi(a)_{(n)}\phi(b)$ for all $a,b\in\A$ and $n\in\Z_+$
\item[{\rm (H2)}] $\p_{\B}\circ \phi = \phi \circ \p_{\A}.$
\end{itemize}
By means of these morphisms we define the category of $\R$-conformal
superalgebras, which we denote by $\R-conf.$

 An $\R$-conformal superalgebra homomorphism
$\phi:\A\rightarrow\B$ is an {\em isomorphism} if it is bijective;
it is an {\em automorphism} if also $\A=\B$.  The set of all
automorphisms of an $\R$-conformal superalgebra $\A$ will be denoted
$\Aut_{\R-conf}(\A)$, or simply $\Aut_{\R}(\A).$

\begin{remark}\label{lambdaproduct}{\rm To simplify some of the longer computations, it will be
convenient to use the $\gl$-product.  This is the generating
function $  a_\lambda  b  $ defined as
\begin{equation*}
  a_\lambda  b  :=\sum_{n\in\Z_+}\gl^{(n)}a_{(n)}b
\end{equation*}
for any pair of conformal superalgebra elements $a,b$ and
indeterminate $\gl$, with $\gl^{(n)}:=\frac{1}{n!}\gl^n$.  In the
case of Lie conformal superalgebras, we denote  $a_\lambda  b$ by $[
a_\lambda  b ].$

The condition (H1) that superconformal homomorphisms $\phi$ respect
all $n$-products is equivalent to
\begin{enumerate}
\item[{\rm (H1')}]$\displaystyle{\phi(a_\lambda b)=\phi(a)_\lambda \phi(b)}$
\end{enumerate}
for all $a,b$.

Axiom (CS0) is equivalent to the property that $a_\lambda  b$ is
polynomial in $\lambda  .$  Axion (CS1) is equivalent to:
$$
\big(\partial  _{\A}(a)\big)_\lambda  b =-\lambda  a_\lambda  b
\q\text{\rm and}\q a_\lambda  \partial  _{\A}b = (\partial
_{\A}+\lambda  )(a_\lambda  b).
$$
Axiom (CS2) is equivalent to
$$
a_\lambda  rb = r(a_\lambda  b)\q\text{\rm and}\q (ra)_\lambda  b =
(a_{\lambda  +\delta  _R}b)_{\to r},
$$
where $\to$ means that $\delta  _R$ is moved to the right and
applied to $r.$}
\end{remark}

\begin{remark}\label{homonot} {\rm Note that
the homotheties $r_\A:\ a\mapsto ra$ (for $r\in R$) are typically
{\em not} $\R$-conformal superalgebra homomorphisms, and that the
map $\p_\A :\ a\mapsto \p_\A (a)$ is a $\R$-conformal superalgebra
derivation of the $\lambda  $-product: $\partial  _{\A}(a_\lambda
b)= (\partial  _{\A}a)_\lambda  b + a_\lambda \partial  _{\A}b.$}
\end{remark}

\subsection{Base change}

Let $\cS = (S,\delta_S)$ be an extension of a base ring $\R = (R,
\delta_R)\in\kalgder$. Given an $\R$-conformal superalgebra $\A$,
the $S$-module $\A\ot_R S$ admits an $\cS$-conformal structure,
which we denote by $\A \ot_\R \cS$, that we now describe.

The derivation $\p_{\A\ot_\R\cS}$ is given by
\begin{equation}\label{extdefi}
\p_{\A \ot_\R \cS}(a\ot s) := \p_{\A}(a)\ot s + a\ot\delta_S(s)
\end{equation}
for all $a\in\A$, $s\in S.$ The $\Z/2Z$-grading is inherited  from
that of $\A$ by setting
$$(\A\ot_{\R}\cS)_{\ol{\iota}}:=\A_{\ol{\iota}}\ot_\R\cS.$$
The $n$-products are defined via
\begin{equation}\label{extproddef}
(a\ot r)_{(n)}(b\ot
s)=\sum_{j\in\Z_+}(a_{(n+j)}b)\ot\delta_S^{(j)}(r)s
\end{equation}
for all $a,b\in\A$, $r,s\in S$, and $n\in\Z_+$. Axioms (CS0)--(CS3)
hold, as can be verified directly. (If $\A$ is also a {\em Lie}
conformal superalgebra, then (CS4)--(CS5) hold, and $\A\ot_\R\cS$ is
also Lie.) The $\cS$-conformal superalgebra on $\A \ot_\R \cS$
described above is said to be obtained from $\A$ by {\it base
change} from $\R$ to $\cS$.

\medskip
\begin{example} {\rm The affinization $\widehat{\A}$ of a $k$-conformal superalgebra $\A$ (Example \ref{loopexample}), viewed as a conformal superalgebra over $\R:=(k[t,t^{-1}],\delta_t)$,
 is obtained from $\A$ by base change from $k$ to $\R$.}
 \end{example}

\medskip

\begin{remark}\label{tensorsareassociative}{\rm It is straightforward to verify that the tensor products
used in defining change of base are associative.  More precisely,
assume that $\cS = (S,\delta_S)$ is an extension of both $\R =
(R,\delta_R)$ and $\cT = (T,\delta_T)$, and $\cU = (U,\delta_U)$ is
also an extension of $\cT = (T,\delta_T)$.  Then for any
$\R$-conformal superalgebra $\A$, the map $(a\ot s)\ot u\mapsto
a\ot(s\ot u)$ defines a $\cU$-conformal isomorphism
$$\left(\A\ot_\R \cS\right)\ot_\cT \cU \cong \A\ot_\R\left(\cS\ot_\cT \cU\right).$$
Here  $u\in U$ acts on $\A\ot_R(S\ot_T U)$ by multiplication, namely
$$u(a\ot(s\ot u')):=a\ot (s\ot uu')$$
for $a\in\A$, $s\in S$ and $u'\in U$; the derivation  $\p_{\A\ot
(\cS\ot \cU)}$ acts on $\A\ot_R(S\ot_T U)$ as
$$\p_{\A\ot_\R (\cS\ot \cU)}=\p_\A\ot \hbox{id}_{S\ot U} + \hbox{id}_\A\ot\delta_{S\ot U},$$
where $\delta_{S\ot U}:=\delta_S\ot \hbox{id}_U +\hbox{id}_S\ot
\delta_U$.  The associativity of tensor products will be useful when
working with $\cS/\R$-forms (\S 2 and \S 3 below).}
\end{remark}

\noindent {\bf Extension functor:} Each extension $\cS =
(S,\delta_S)$ of $\R = (R,\delta_R)$ defines an {\it extension
functor}
$$ \fE = \fE_{\cS/\R} : \R-conf \to \cS-conf$$
as follows:

Given an $\R$-conformal superalgebra $\A$, let $\fE(\A)$ be the
$\cS$-conformal superalgebra $\A\ot_\R \cS.$ For each $\R$-conformal
superalgebra homomorphism $\psi:\ \A\rightarrow {\mathcal B}$, the
unique $S$-linear map satisfying
\begin{eqnarray}
\fE(\psi):\ \fE(\A)&\rightarrow& \fE({\mathcal B})\\
a\ot s&\mapsto&\psi(a)\ot s
\end{eqnarray}
is clearly a homomorphism of $\cS$-conformal superalgebras, and it
is straightforward to verify that $\fE$ is a functor.

\bigskip

\noindent {\bf Restriction functor:} Likewise, any $\cS$-conformal
superalgebra $\B$ can be viewed as an $\R$-conformal superalgebra by
{\it restriction of scalars from} $\cS$ to $\R$:

If the extension $\cS/\R$ corresponds to a $\kalgder$ morphism
$\phi:\ R\rightarrow S$, we view $\B$ as an $R$-module via $\phi$.
Then $\B$ is naturally an $\R$-conformal superalgebra. The only
nontrivial
 axiom to verify is (CS2). Using the notation of Remark \ref{homonot}, we have
\begin{eqnarray*}
\p_{\B}\circ r_{\B} &=& \p_{\B}\circ\phi(r)_{\B}\\
&=& \phi(r)_{\B}\circ\p_{\B}+\delta_S(\phi(r))_{\B}\\
&=& r_{\B}\circ\p_{\B}+ \phi
(\delta_R(r))_{\B}\\
 &=& r_{\B}\circ\p_{\B}+ \delta_R(r)_{\B}.
 \end{eqnarray*}

This leads to the {\it restriction functor}
$$ \fR = \fR_{\cS/\R} : \cS-conf \to \R-conf,$$
which attaches to an $\cS$-conformal superalgebra $\B$ the same $\B$
viewed as an $\R$-conformal superalgebra; likewise, to any
$\cS$-superconformal homomorphism $\psi:\ \B \rightarrow\cC$, $\fR$
attaches the $\R$-conformal superalgebra morphism $\psi$.

\subsection{The automorphism functor of a conformal superalgebra}

Let $\A$ be an $\R = (R, \delta_R)$-conformal superalgebra. We now
define the automorphism group functor ${\bf Aut}(\A)$. For each
extension $\cS = (S, \delta_S)$ of $\R$, consider the group
\begin{equation}
{\bf Aut}(\A)(\cS):=\Aut_{\cS}(\A\ot_\R \cS)
\end{equation}
of automorphisms of the $\cS$-conformal superalgebra $\A \ot_\R
\cS$. For each morphism $\psi: \cS_1\rightarrow \cS_2$ between two
extensions $\cS_1 = (S_1 , \delta_1)$ and $\cS_2 = (S_2 ,\delta_2)$
of $\R$, and each automorphism $\theta\in{\bf Aut}(\A)(\cS_1)$, let
${\bf Aut}(\A)(\psi)(\theta)$ be the unique $S_2$-linear map
determined by
\begin{eqnarray}
{\bf Aut}(\A)(\psi)(\theta):\ \A\ot_R S_2&\rightarrow&\A\ot_R S_2\\
a\ot 1&\mapsto&\sum_ia_i\ot\psi(s_i)
\end{eqnarray}
for $a\in\A$, where $\theta(a\ot 1)=\sum_ia_i\ot s_i.$

\begin{proposition}\label{autfunctor} ${\bf Aut}(\A)$ is a functor from the category  of extensions of $(R,\delta_R)$ to
the category of groups.
\end{proposition}
\proof Let $\theta_1,\theta_2\in\Aut_{\cS_1}(\A\ot_\R \cS_1)$, and
write $\theta_2(a\ot1)=\sum_ia_i\ot s_i$ for some $a_i\in\A$ and
$s_i\in S_1$.  Then for any morphism $\psi : \cS_1 \to \cS_2$, we
have (in the notation above):

\smallskip
\noindent ${\bf Aut}(\A)(\psi)(\theta_1)\circ{\bf
Aut}(\A)(\psi)(\theta_2)(a\ot 1)$
\begin{eqnarray*}
\hspace{1in}&=&{\bf Aut}(\A)(\psi)(\theta_1)(1\ot\psi)\theta_2(a\ot 1)\\
&=&{\bf Aut}(\A)(\psi)(\theta_1)\sum_ia_i\ot\psi(s_i)\\
&=&\sum_i\psi(s_i){\bf Aut}(\A)(\psi)(\theta_1)(a_i\ot 1)\\
&=&(1\ot\psi)\sum_is_i\theta_1(a_i\ot 1)\\
&=&(1\ot\psi)\theta_1\sum_ia_i\ot s_i\\
&=&(1\ot\psi)\theta_1\theta_2(a\ot 1)\\
&=&{\bf Aut}(\A)(\psi)(\theta_1\theta_2)(a\ot 1).
\end{eqnarray*}
Using the $S_2$-linearity of the $\cS_2$-conformal automorphisms
${\bf Aut}(\A)(\psi)(\theta_1)$, ${\bf Aut}(\A)(\psi)(\theta_2)$,
and ${\bf Aut}(\A)(\psi)(\theta_1\theta_2)$, we have
$${\bf
Aut}(\A)(\psi)(\theta_1)\circ{\bf Aut}(\A)(\psi)(\theta_2)={\bf
Aut}(\A)(\psi)(\theta_1\theta_2).$$

In particular, note that for any $\theta\in{\bf Aut}(\A)(\cS_1)$, we
have
\begin{equation}
{\bf Aut}(\A)(\psi)(\theta^{-1})\circ{\bf
Aut}(\A)(\psi)(\theta)={\bf Aut}(\A)(\psi)({\rm id}_{\A\ot_\R
\cS_1}).
\end{equation}
 It is clear from
the definition of ${\bf Aut}(\A)(\psi)$ that ${\bf
Aut}\,\A(\psi)({\rm id}_{\A\ot_\R \cS_1})$ is the identity map on
$\A\ot 1$, hence also on $\A\ot_R S_2$ by $S_2$-linearity.  Thus
${\bf Aut}(\A)(\psi)(\theta)$ has a left inverse and is therefore
injective.  Interchanging the roles of $\theta$ and $\theta^{-1}$
shows that ${\bf Aut}(\A)(\psi)(\theta)$ has a right inverse, so it
is also surjective.

That ${\bf Aut}(\A)(\psi)(\theta)$ is an $\cS_2$-conformal
superalgebra homomorphism follows easily from the assumption that
$\theta\in\Aut_{\cS_1}(\A\ot_\R \cS_1)$ and $\psi: \cS_1 \rightarrow
\cS_2$ is a morphism of $\R$-extensions.  Therefore, ${\bf
Aut}(\A)(\psi)(\theta)\in\Aut_{\cS_2}(\A\ot_\R \cS_2)$, and we have
now shown that ${\bf Aut}(\A)(\psi)$ is a group homomorphism
\begin{equation}\
{\bf Aut}(\A)(\psi):\ {\bf Aut}(\A)(\cS_1)\rightarrow {\bf
Aut}(\A)(\cS_2).
\end{equation}

Clearly  ${\bf Aut}(\A)$ sends the identity morphism ${\rm
id}_{\cS}$ to the identity map on ${\bf Aut}(\A)(\cS)$ for any
extension $\cS$ of $\R$. To finish proving that ${\bf Aut}(\A)$ is a
functor, it remains only to note that if $\psi_1:\ \cS_1 \rightarrow
\cS_2 $ and $\psi_2:\ \cS_2 \rightarrow \cS_3$ are morphisms between
extensions $\cS_i = (S_i, \delta_i)_{1 \leq i \leq 3}$ of $\R$, then
for $a\in\A$ and $\theta(a\ot 1)=\sum_ia_i\ot s_i$, we have
$${\bf Aut}(\A)(\psi_2\psi_1)(\theta):\ a\ot 1\mapsto\sum_ia_i\ot\psi_2\psi_1(s_i),$$
which defines precisely the same map (via $S_3$-linearity) as ${\bf
Aut}(\A)(\psi_2)\circ{\bf Aut}(\A)(\psi_1)(\theta).$  Hence ${\bf
Aut}(\A)(\psi_2\psi_1)={\bf Aut}(\A)(\psi_2)\circ{\bf
Aut}(\A)(\psi_1)$, which completes the proof of the proposition.\qed

\bigskip

\section{Forms of conformal superalgebras and \v Cech cohomology}\label{ketto}

 Given $R$ in $\kalg$ and a (not necessarily commutative, associative, or unital)
$R$--algebra $A$, recall that a {\em form of} $A$ (for the
$\fppf$--topology on $\Spec(R)$) is an $R$-algebra $F$ such that
$F\ot_R S \cong A\ot_R S$ (as $S$-algebras) for some fppf
(faithfully flat and finitely presented) extension $S/R$ in $\kalg.$
There is a correspondence between $R$-isomorphism classes of forms
of $A$ and the pointed set of non-abelian cohomology $H_{\fppf}^1(R,
\bAut(A))$ defined \`a la \v Cech. Here $\bAut(A):=\bAut(A)_R$
denotes the sheaf of groups over $\Spec(R)$ that attaches to an
extension $R'/R$ in $\kalg$ the group $\Aut_S(A\ot_R R')$ of
automorphisms of the $R'$--algebra $A \otimes_R R'.$ For any
extension $S/R$ in $\kalg$, there is a canonical map

$$ H_{\fppf}^1(R, \bAut(A)) \to H_{\fppf}^1(S, \bAut(A)_S)$$
The kernel of this map is denoted by $H_{\fppf}^1(S/R, \bAut(A));$
these are the forms of $A$ that are trivialized by the base change
$S/R.$ One has

$$ H_{\fppf}^1(R, \bAut(A))  = \limind H_{\fppf}^1(S/R, \bAut(A)),$$
where the limit is taken over all fppf extensions $S/R$ in $\kalg$.

In trying to recreate this construction for an $\R$-conformal
superalgebra $\A$, we encounter a fundamental obstacle: Unlike in
the case of algebras, the $n$-products (\ref{extproddef}) in
$\A\ot_\R \cS$ are {\it not} obtained by $S$-linear extension of the
$n$-products in $\A$ (unless the derivation of $\cS$ is trivial).
This prevents the automorphism functor ${\bf Aut}(\A)$ from being
representable in the na\"ive way, and the classical theory of forms
cannot be applied blindly. Nonetheless, we will show in the next
section that the expected correspondence between forms and
cohomology continues to hold, even in the case of conformal
superalgebras.

In the case of algebras, when the extension $S/R$ is Galois,
isomorphism classes of $S/R$--forms have an interpretation in terms
of non-abelian Galois cohomology. (See \cite{waterhouse}, for
instance.)
We will show in \S \ref{negy} that, just as in the case of algebras,
the Galois cohomology $H^1\left(\Gal(\cS/\R), \Aut_\cS(\A\ot_\R
\cS)\right)$ still parametrizes the $\cS/\R$-forms of $\A$ (with the
appropriate definition of Galois extension and $\Aut_\cS(\A\ot_\R
\cS)$).

 Throughout this section, $\A$ will denote a conformal superalgebra over $\R = (R, \delta_R)$.

\subsection{Forms split by an extension}

\begin{definition}{\em Let  $\cS $ be an extension of $\R$.  An $\R$-conformal superalgebra $\F$ is an
{\em $\cS/\R$-form} of $\A$ (or {\em form of $\A$ split by $\cS$})
if
$$\F\ot_\R \cS\cong\A\ot_\R \cS$$
as $\cS$-conformal superalgebras. }\end{definition}

For us, the most interesting examples of forms split by a given
extension are the conformal superalgebras that are obtained via the
type of twisted loop construction that one encounters in the theory
of affine Kac-Moody Lie algebras.

\begin{example}\label{loopalgdef}{\em Assume $k$ is algebraically closed.
Suppose that $\A$ is a $\K$-conformal superalgebra, equipped with an
automorphism $\gs$ of period $m$. For each $i \in \Z$  consider the
eigenspace
$$\A_{i}=\{x\in\A\ |\ \gs(x)=\xi_m^ix\}.$$
with respect to our fixed choice $(\xi_m)$ of compatible primitive
roots of unity in $k.$ (The space $\A_i$ depends only on the class
of $i$ modulo $m$, of course.) Let $R=\K[t,t^{-1}]$ and
$S_m=\K[t^{1/m},t^{-1/m}]$, and let $\cS_m = (S_m, \delta_t)$ be the
extension of $\R = (R, \delta_t)$ where $\delta_t = \frac{d}{dt}.$

Consider the subspace $\Loop(\A,\gs) \subseteq \A \ot_k \cS_m$ given
by
\begin{equation}
\Loop(\A,\gs)=\bigoplus_{i\in\Z}\A_{i}\ot t^{i/m}.
\end{equation}
 Each eigenspace $\A_i$ is stable under $\p_\A$ because $\sigma$
is a conformal automorphism. From this, it easily follows that
$\Loop(\A,\gs)$ is stable under the action of $\p_{\A \ot \cS_m} =
\p_\A \ot 1 + 1 \ot \delta_t.$ Since $\Loop(\A,\gs)$ is also closed
under the $n$-products of $\A \ot_k \cS_m$ , it is a $k$-conformal
subalgebra of $\A \ot_k \cS_m$ called  the {\it (twisted) loop
algebra} of $\A$ with respect to $\gs.$ (Note that the definition of
$\Loop(\A,\gs)$ does not depend on the choice of the period $m$ of
the given automorphism $\sigma$).
 It is clear that $\Loop(\A,\gs)$ is stable under the natural action of
 $R$ on $\A \ot_k \cS_m.$
 As in Example \ref{loopexample}, one checks that $\p_{\A \ot \cS_m} \circ r_{\A \ot_k \cS_m} - r_{\A \ot_k
\cS_m} \circ \p_{\A \ot \cS_m} = \delta_t(r)_{\A \ot_k \cS_m}$ for
all $r \in R$. This shows that just as in the untwisted case, the
twisted loop algebra $\Loop(\A,\gs)$ is both a $k$- and an
$\R$-conformal superalgebra.}

\end{example}

\begin{proposition}\label{loopsareforms}  Let  $\gs$ be an automorphism of period $m$ of a $\K$-conformal superalgebra $\A$.
Then the twisted loop algebra $\Loop(\A,\gs)$ is an $\cS_m/\R$-form
of $\A\ot_\K \R$.
\end{proposition}
\proof For ease of notation, we will write $\cS=(S,\delta_S)$ for
$\cS_m=(S_m,\delta_{S_m})$ in this proof.  By the associativity of
the tensor products used in scalar extension (Remark
\ref{tensorsareassociative}), the multiplication map
\begin{eqnarray}
\psi:\ (\A\ot_\K \R)\ot_\R \cS&\rightarrow& \A\ot_\K \cS\\
(a\ot r)\ot s&\mapsto& a\ot rs
\end{eqnarray}
(for $a\in\A$, $r\in R$, and $s\in S$) is an isomorphism of
$\cS$-conformal superalgebras.

Likewise, it is straightforward to verify that the multiplication
map
\begin{eqnarray}
\mu:\ \big(\A\ot_\K \cS\big)\ot_\R \cS&\rightarrow& \A\ot_\K \cS\\
(a\ot s_1)\ot s_2&\mapsto&a\ot s_1s_2
\end{eqnarray}
(for $a\in\A$ and $s_1,s_2\in S$) is a homomorphism of
$\cS$-conformal superalgebras.  Indeed, $\mu$ is the composition of
the ``associativity isomorphism''
$$\big(\A\ot_\K \cS\big)\ot_\R \cS\rightarrow \A\ot_\K\big(\cS\ot_\R \cS)$$
with the superconformal homomorphism defined by multiplication:
\begin{eqnarray*}
\A\ot_\K(\cS\ot_\R \cS)&\rightarrow& \A\ot_\K \cS\\
a\ot (s_1\ot s_2)&\mapsto& a\ot s_1s_2.
\end{eqnarray*}

To complete the proof of Proposition \ref{loopsareforms}, it
suffices to prove that the restriction
\begin{equation}
\mu:\ \Loop(\A,\gs) \ot_\R \cS\rightarrow \A\ot_\K \cS
\end{equation}
is bijective. For $a_i \in\A_i$, we have
$$a_i \ot t^{j/m}=\mu(a_i \ot t^{i/m}\ot t^{(j-i)/m}),$$
so $\mu$ is clearly surjective. To see that $\mu$ is also injective,
assume (without loss of generality) that a $\K$-basis $\{a_\gl\}$ of
$\A$ is chosen so that each $a_{\gl} \in \A_{i(\lambda)}$ for some
unique $0 \leq i(\lambda) < m.$ Let $x\in \Loop(\A,\gs) \ot_R S$.
Since $S$ is a free $R$-module with basis $\{t^{i/m}\ |\ 0\leq i\leq
m-1\}$, we can uniquely write
$$x=\sum_{i=0}^{m-1}x_i\ot t^{i/m}$$
where $x_i=\sum_\gl a_\gl\ot f_{\gl i}$ and $f_{\gl i}\in
t^{i(\gl)/m}\K[t,t^{-1}]$.  Then if $\mu(x)=0$, we have
$$\sum_{\gl i}a_\gl\ot f_{\gl i}t^{i/m}=0,$$
and thus
$$\sum_{i=0}^{m-1} f_{\gl, i}t^{i/m}=0$$
for all $\gl$.  Then $f_{\gl, i}=0$ for all $\gl$ and $i$.   Hence
$x=0$, so $\mu$ is injective, and
\begin{equation}
\psi^{-1}\circ\mu:\ \Loop(\A,\gs) \ot_\R \cS\rightarrow\big(\A\ot_\K
\R\big)\ot_\R \cS
\end{equation}
is an $\cS$-conformal superalgebra isomorphism as desired. \qed

\subsection{Cohomology and forms}\label{negy}

Throughout this section $\cS = (S,\delta_S)$ will denote  an
extension of $\R = (R,\delta_R)$, and $\A$ will be an $\R$-conformal
superalgebra.

\begin{lemma}\label{preservesnprod} Let $\psi:\ \A\rightarrow\B$ be an $\R$-conformal superalgebra
homomorphism and $\gamma:\ \cS\rightarrow \cS$  a morphism of
extensions. The canonical map
\begin{equation}
\psi\ot\gamma:\ \A\ot_R S\rightarrow \B\ot_R S
\end{equation}
is $R$-linear, commutes with the action of $\p_{\A\ot_\R \cS},$ and
preserves $n$-products. In particular, $\psi\ot\gamma$ is an
$\R$-conformal superalgebra homomorphism via restriction of scalars
from $\cS$ to $\R$:
$$\psi\ot\gamma:\ \fR_{\cS/\R}(\A\ot_\R \cS)\rightarrow\fR_{\cS/\R}(\B\ot_\R \cS).$$
\end{lemma}
\proof Let $x\in\A$ and $s\in S$. Then
\begin{eqnarray*}
\psi\ot\gamma\big(\p_{\A \ot_\R \cS}(x\ot s)\big)&=&\psi\ot\gamma\big(\gd_{\A}(x)\ot s+ x\ot\delta_S(s)\big)\\
&=&\p_\B(\psi(x))\ot\gamma(s)+\psi(x)\ot\delta_S\big(\gamma(s)\big)\\
&=&\p_{\B \ot_\R \cS}\big(\psi(x)\ot\gamma(s)\big).
\end{eqnarray*}
Also, for $x,y\in\A$ and $s,t\in S$, we have
\begin{eqnarray*}
\psi\ot\gamma (x \ot s_{(n)}y\ot t)&=&\psi\ot\gamma\left(\sum_{j\in\Z_+}x_{(n+j)}y\ot\delta_S^{(j)}(s)t\right)\\
&=&\sum_{j\in\Z_+}\psi(x)_{(n+j)}\psi(y)\ot\delta_S^{(j)}\big(\gamma(s)\big)\gamma(t)\\
&=&\psi(x)\ot\gamma(s)_{(n)}\psi(y)\ot\gamma(t).
\end{eqnarray*}\qed

\begin{corollary}\label{psione}
The map $\psi\ot 1:\ \A\ot_\R \cS\rightarrow \B\ot_\R \cS$ is an
$S$-conformal superalgebra homomorphism.
\end{corollary}
\proof It is enough to note that the map $\psi\ot 1$ commutes with
the action of $S.$ \qed
\bigskip

  For $1\leq i\leq 2$ and $1\leq j< k\leq
3$, consider the following $R$-linear maps:
\begin{eqnarray*}
&&d_i:\ S\rightarrow S\ot_R S\\
&&d_{jk}:\ S\ot_R S\rightarrow S\ot_R S\ot_R S,
\end{eqnarray*}
defined by $d_1(s)=s\ot 1$, $d_2=1\ot s$, $d_{12}(s\ot t)=s\ot t\ot
1$, $d_{13}(s\ot t)=s\ot 1\ot t$, and $d_{23}(s\ot t)=1\ot s\ot t$
for all $s,t\in S$.  It is straightforward to verify that these
induce $\R-ext$ morphisms $d_i:\ \cS\rightarrow \cS\ot_\R \cS$ and
$d_{jk}:\ \cS\ot_\R \cS\rightarrow \cS\ot_\R \cS\ot_\R \cS$.  (See
\S 1.1). Let $\A$ be an $\R$-conformal superalgebra. By
functoriality (Proposition \ref{autfunctor}), we obtain group
homomorphisms (also denoted by $d_i$ and $d_{jk}$)
\begin{eqnarray*}
&&d_i:\ {\bf Aut}(\A) (\cS)\rightarrow {\bf Aut}(\A)(\cS\ot_\R \cS)\\
&&d_{jk}:\ {\bf Aut}(\A)(\cS\ot_\R \cS)\rightarrow {\bf
Aut}(\A)(\cS\ot_\R \cS\ot_\R \cS)
\end{eqnarray*}
for $1\leq i\leq 2$ and $1\leq j<k\leq 3$.

Recall that $u\in {\bf Aut}(\A)(\cS\ot_\R \cS)$ is called a {\em
$1$-cocycle}\footnote{For faithfully flat ring extensions $S/R$, the
$1$-cocycle condition is motivated by patching data on open
coverings of $\Spec(R)$.
See the discussion in \cite[\S 17.4]{waterhouse}, for instance.} if
\begin{equation}\label{cocyclecondition}
d_{13}(u)=d_{23}(u)d_{12}(u).
\end{equation}
On the set $Z^1(\cS/\R, {\bf Aut}(\A))$ of $1$-cocycles, one defines
an equivalence relation by declaring two cocycles $u$ and $v$ to be
{\em equivalent} (or {\em cohomologous}) if there exists an
automorphism $\gl\in {\bf Aut}(\A)(\cS)$ such that
\begin{equation}
v=\big(d_2(\gl)\big)\,u\,\big(d_1(\gl)\big)^{-1}.
\end{equation}
The corresponding quotient set is denoted $H^1(\cS/\R, {\bf
Aut}(\A))$ and is the (nonabelian) \v Cech cohomology relative to
the covering $\Spec(S)\rightarrow \Spec(R)$.  There is no natural
group structure on this set, but it has a distinguished element,
namely the equivalence class of the identity element of the group
${\bf Aut}(\A)(\cS\ot_\R \cS)$.   We will denote this class by $1$,
and write $[u]$ for the equivalence class of an arbitrary cocycle
$u\in Z^1(S/R, \autfun(\A))$.

\begin{theorem}\label{forms}  Assume that the extension $\cS = (S,\delta_S)$ of $\R = (R,\delta_R)$ is
 faithfully flat (i.e.,~$S$ is a faithfully flat $R$-module).  Then for any $\R$-conformal
superalgebra $\A$, the pointed set $H^1(\cS/\R, {\bf Aut}(\A))$
parametrizes the set of $\R$-isomorphism classes of $\cS/\R$-forms
of $\A$.  Under this correspondence, the distinguished element $1$
corresponds to the isomorphism class of $\A$ itself.
\end{theorem}
\proof It suffices to check that the standard descent formalism for
modules is compatible with the conformal superalgebra structures.

Throughout this proof, fix an $\cS/\R$-form $\B$ of the
$\R$-conformal superalgebra $\A$.  Let $\eta:\ S\ot_R S\rightarrow
S\ot_R S$ be the ``switch'' map given by $\eta(s\ot t)=t\ot s$ for
all $s,t\in S$. Let
\begin{eqnarray*}
\eta_\A:=\hbox{id}_\A\ot\eta:\ \A\ot_R S\ot_R S&\rightarrow& \A\ot_R S\ot_R S\\
\eta_\B:=\hbox{id}_\B\ot\eta:\ \B\ot_R S\ot_R S&\rightarrow& \B\ot_R
S\ot_R S.
\end{eqnarray*}

\bigskip

We now check that the key points in the classical faithfully flat
descent formalism hold in the conformal setting:

\begin{enumerate}
\item[{\rm \bf (1)}] {\it Let $\psi:\ \B\ot_\R \cS\rightarrow
\A\ot_\R \cS$ be an isomorphism of $\cS$-conformal superalgebras.
Let
\begin{equation}
u_{\psi,\B}:=(\eta_\A)(\psi\ot 1)(\eta_\B)(\psi^{-1}\ot 1).
\end{equation}
Then $u_{\psi,\B}\in Z^1(\cS/\R, {\bf Aut}(\A))$.}

\end{enumerate}
\noindent {\em Proof (1):}  That $u_{\psi,\B}:\ \A\ot_R S\ot_R
S\rightarrow \A\ot_R S\ot_R S$ is $S\ot_R S$-linear and bijective is
 clear, and it is straightforward to verify that $u_{\psi,\B}$
satisfies the cocycle condition (\ref{cocyclecondition}).  Therefore
it is enough to check that $u_{\psi,\B}$ commutes with the
derivation $\p_{\A \ot_\R \cS\ot_R \cS}$, and that it preserves the
$n$-products.

By Corollary \ref{psione} and the asociativity of the tensor
product, $\psi\ot 1$ and $\psi^{-1}\ot 1$ are $\cS \otimes_R
\cS$-superconformal homomorphisms, so it is only necessary to check
that $\eta_\A$ and $\eta_\B$ commute with $\p_{\A \ot_\R \cS\ot_R
\cS}$ and preserve $n$-products. By Lemma \ref{preservesnprod},
applied to $\eta_\A=\hbox{id}_\A\ot\eta$, it is enough to check that
$\eta$ commutes with $\delta_{S\ot S},$ which is clear since

\begin{eqnarray}\label{etacommutes}
\eta(\delta_{S\ot S}(s\ot t))&=&\eta\big(\delta_S(s)\ot t+s\ot \delta_S(t)\big)\\
&=&t\ot \delta_S(s)+ \delta_S(t)\ot s\label{etacommutes1}\\
&=&\delta_{S\ot S}\big(\eta_\A(s\ot t)\big)\label{etacommutes2}
\end{eqnarray}
for all $s,t\in S$.

\bigskip
\bigskip

\begin{enumerate}

\item[{\rm \bf (2)}] {\it The class of $u_{\psi,\B}$ in
$H^1(\cS/\R,\autfun(\A))$ is independent of the choice of
automorphism $\psi$ in Part (1).  If we denote this class by
$[u_\B]$, then $\gs:\ \B\mapsto [u_\B]$ is a map from the set of
$\R$-isomorphism classes $\cS/\R$-forms of $\A$ to the pointed set
$H^1(\cS/\R,\autfun(\A))$.}
\end{enumerate}
\noindent {\em Proof (2):}  Suppose that $\phi$ is another
$\cS$-superconformal isomorphism $$\phi:\ \B\ot_\R \cS\rightarrow
\A\ot_\R \cS.$$ Let $\gl=\phi\psi^{-1}\in\autfun(\A)(S)$.  Note that
$d_2(\gl)\eta_\A=\eta_\A d_1(\gl)$.  Thus
\begin{eqnarray*}
d_2(\gl)u_{\psi,\B}d_1(\gl)^{-1}&=&d_2(\gl)\eta_\A(\psi\ot 1)\eta_\B(\psi^{-1}\ot 1)(\psi\phi^{-1}\ot 1)\\
&=&\eta_\A(\phi\psi^{-1}\ot 1)(\psi\ot 1)\eta_\B(\phi^{-1}\ot 1)\\
&=&u_{\phi,\B}.
\end{eqnarray*}

\bigskip

\begin{enumerate}
\item[{\rm \bf (3)}] {\it Let $u\in\autfun(\A)(\cS\ot_\R \cS)$.
Then
\begin{enumerate}
\item[{\rm \bf (a)}] The subset $\A_u:=\{x\in\A\ot_R S\ |\ \ u(x\ot 1)=\eta_\A(x\ot 1)\}$ is an
$\R$-conformal subalgebra of $\A\ot_\R \cS$.
\item[{\rm \bf (b)}] The canonical map
\begin{eqnarray*}
\mu_u:\ \A_u\ot_\R \cS&\rightarrow& \A\ot_\R \cS\\
x\ot s&\mapsto&s.x
\end{eqnarray*}
is an $\cS$-conformal superalgebra isomorphism.
\item[{\rm \bf (c)}] If $u$ and $v$ are cohomologous cocycles in
$Z^1(\cS/\R,\autfun(\A))$, then $\A_u$ and $\A_v$ are isomorphic as
$\R$-conformal superalgebras.
\end{enumerate}}
\end{enumerate}
\noindent {\em Proof (3a):}  Clearly $\A_u$ is an $R$-submodule of
$\A\ot_R S$.  Next we verify that $\A_u$ is stable under the action
of $\p_{\A \ot_\R \cS}$.

Recall (\ref{etacommutes}) that  $\eta_\A$  commutes with the
derivation $ \p_{\A \ot_\R \cS \ot_\R  \cS}.$ Thus for all
$x\in\A_u$,
\begin{eqnarray*}
 u(\p_{\A \ot_\R \cS}(x)\ot1)&=& u\,\p_{\A \ot_\R \cS \ot_\R  \cS}(x\ot 1)\\
&=&\p_{\A \ot_\R \cS \ot_\R  \cS}\,u(x\ot 1)\\
&=&\p_{\A \ot_\R \cS \ot_\R  \cS}\,\eta_\A(x\ot 1)\\
&=&\eta_\A \p_{\A \ot_\R \cS \ot_\R  \cS}(x\ot 1)\\
&=&\eta_\A \p_{\A \ot_\R \cS}(x)\ot 1.
\end{eqnarray*}
To complete the proof of (3a), it remains only to show that $\A_u$
is closed under $n$-products. For $x$ and $y$ in $\A_u$ we have
\begin{eqnarray*}
u\big((x_{(n)}y)\ot 1\big)&=&u(x\ot 1_{(n)}y\ot 1)\\
&=&u(x\ot 1)_{(n)}u(y\ot 1)\\
&=&\eta_\A(x\ot 1)_{(n)}\eta_\A(y\ot 1)\\
&=&\eta_\A\big(x\ot 1_{(n)}y\ot 1\big)\\
&=&\eta_\A\big((x_{(n)}y)\ot 1\big)
\end{eqnarray*}
for all $x,y\in\A_u$ \big(where we have used (\ref{etacommutes}),
(\ref{etacommutes1}), and (\ref{etacommutes2}) to get
\newline $\eta_\A(x\ot 1)_{(n)}\eta_\A(y\ot 1)=\eta_\A\big(x\ot
1_{(n)}y\ot 1\big)$\big).

\bigskip

\noindent {\em Proof (3b):} The map $\mu_u:\ \A_u\ot_R
S\rightarrow\A\ot_R S$ is an $S$-module isomorphism by the classical
descent theory for modules.  (See \cite[Chap 17]{waterhouse}, for
instance.)  We need only show that it is a homomorphism of
$\cS$-conformal superalgebras.

Let $\mu$ be the multiplication map
\begin{eqnarray*}
\mu:\ \A\ot_R S\ot_R S&\rightarrow& \A\ot_R S\\
a\ot s\ot t&\mapsto& a\ot st.
\end{eqnarray*}
It is straightforward to verify that
\begin{equation}\label{mu}
\mu\circ\p_{\A \ot_\R \cS \ot \cS}=\p_{\A \ot_\R \cS}\circ\mu, \,
{\rm and} \, \mu \, {\rm preserves} \, n{\rm-products.}
\end{equation}
 Since $\mu_u$ is the restriction of $\mu$ to
$\A_u\ot_R S$ it follows from (\ref{mu}) that $\mu_u$ preserves
$n$-products. It remains only to show that $\mu_u$ commutes with
$\p_{\A_u \ot_\R \cS}$. But since $\p_{\A_u \ot_\R \cS}$ acts on
$\A_u\ot_R S$ as the restriction of the derivation $\p_{\A \ot_\R
\cS \ot_\R \cS}$ of $\A\ot_\R \cS \ot_\R \cS$ to the subalgebra
$\A_u\ot_\R \cS$ and $\mu_u=\mu\,\iota$, where $\iota$ is the
inclusion map $\iota:\ \A_u\ot_R S\hookrightarrow\A\ot_R S\ot_R S$,
we can again appeal to (\ref{mu}). This shows that $\mu_u$ is an
$\cS$-conformal superalgebra isomorphism.



\bigskip

\noindent {\em Proof (3c)}
Write $v=d_2(\gl)\,u\,d_1(\gl)^{-1}$ for some
$\gl\in\autfun(\A)(\cS)$. Then since
$d_2(\gl)\,\eta_\A=\eta_\A\,d_1(\gl)$, we see that
\begin{eqnarray*}
v\,(\gl\ot 1)(a_u\ot 1)&=&v\,d_1(\gl)(a_u\ot 1)\\
&=&d_2(\gl)\,u(a_u\ot 1)\\
&=&d_2(\gl)\,\eta_\A(a_u\ot 1)\\
&=&\eta_\A\,d_1(\gl)(a_u\ot 1)\\
&=&\eta_\A\,(\gl\ot 1)(a_u\ot 1)
\end{eqnarray*}
for all $a_u\in\A_u$.  Thus,
\begin{equation}\label{eqn499}
\gl(\A_u)\subseteq \A_v.
\end{equation}
Likewise, $\gl^{-1}(\A_v)\subseteq\A_u$, so applying $\gl^{-1}$ to
both sides of (\ref{eqn499}) gives:
$$\A_u\subseteq\gl^{-1}(\A_v)\subseteq\A_u,$$
and $\gl(\A_u)=\A_v$.

Since $\gl$ is an $\cS$-automorphism of $\A\ot_\R \cS$, it commutes
with the actions of $S$ and $\p_{\A \ot_\R \cS}$, and it preserves
$n$-products. Thus its restriction to $\A_u$ commutes with $\R$ and
preserves $n$-products. Furthermore $\gl \p_{\A_u} = \p_{\A_v} \gl$
since $\p_{\A_u}$ and $\p_{\A_u}$ are the restrictions of $\p_{\A
\ot_\R \cS}$ to $\A_u$ and $\A_v$ respectively. Thus $\gl$ is an
$\R$-conformal superalgebra isomorphism from $\A_u$ to $\A_v$.

\bigskip

By (3c), there is a well-defined map $\beta:\ [u]\mapsto[\A_u]$ from
the cohomology set $H^1(\cS/\R,\autfun(\A))$ to the set of
$\R$-isomorphism classes of $\cS/\R$-forms of $\A$.  We have also
seen that $\beta$ and the map $\ga$ defined in (2) are inverses of
each other.  That the distinguished element $1\in
H^1(\cS/\R,\autfun(\A))$ corresponds to the algebra $\A$ is clear
from the definition of $\A_u$ given in (3a). \qed

\begin{remark}  {\rm For any isomorphism $\psi$ as in Part (1) of the proof of Theorem \ref{forms},
the cocycle $u_{\psi,\B}$ can be rewritten in terms of the maps
$d_i:\ S\rightarrow S\ot_R S$:}
\begin{equation}
u_{\psi,\B}=d_2(\psi)d_1(\psi)^{-1}.
\end{equation}
\end{remark}
\proof  If $\psi(b\ot 1)=\sum_ia_i\ot w_i$, then
\begin{eqnarray*}
\eta_\A(\psi\ot 1)\eta_\B(b\ot s\ot t)&=&\eta_\A(\psi\ot 1)(b\ot t\ot s)\\
&=&\eta_\A\big(\psi(b\ot t)\ot s\big)\\
&=&\eta_\A\big((t.\psi(b\ot 1))\ot s\big)\\
&=&\eta_\A\left(\sum_ia_i\ot tw_i\ot s\right)\\
&=&\sum_ia_i\ot s\ot tw_i\\
&=&d_2(\psi)(b\ot s\ot t)
\end{eqnarray*}
for all $s,t\in S$.  Thus
\begin{eqnarray*}
u_{\psi,\B}&=&d_2(\psi)(\psi^{-1}\ot 1)\\
&=&d_2(\psi)(\psi\ot 1)^{-1}\\
&=&d_2(\psi)d_1(\psi)^{-1}.
\end{eqnarray*}
\qed

\begin{remark}
{\rm In the notation of Part (1) of the proof of Theorem
\ref{forms}, the algebra $\B$ is isomorphic to the $\R$-conformal
superalgebra $\B\ot 1\subseteq \B\ot_\R \cS$ by the faithful
flatness of $S/R$. This means that there is an isomorphic copy of
each $\cS/\R$-form of $\A$ {\em inside} $\A\ot_\R \cS$, and the
algebra $\B$ can be recovered (up to $\R$-conformal isomorphism)
from the cocycle $u_{\psi,\B}$, since}
\begin{equation}
\psi(\B\ot 1)=\{x\in\A\ot_R S\ |\ u_{\psi,\B}(x\ot 1)=\eta_\A(x\ot
1)\}.
\end{equation}
\end{remark}

\begin{remark} \label{3.17}{\rm Let $\cS = (S,\delta_S)$ be an extension of $\R =
(R,\delta_R).$  Let $\Gamma  $  be a finite group of automorphisms
of $\cS$ (as an extension of $\R)$.  We say that $\cS$ is a {\it
Galois extension of $\R$ with group $\Gamma  $} if $S$ is a Galois
extension of $R$ with group $\Gamma  .$ (See \cite{KO} for
definition). The unique $R$-module map
$$
\psi  :S\otimes_R\,S\to S\times \cdots \times S \quad\quad (\vert
\Gamma  \vert  \text{\rm \ copies of $S$})
$$
satisfying
$$
a\otimes b\mapsto \big(\gamma  (a)b\big)_{\gamma  \in \Gamma  }
$$
is easily seen to be an isomorphism of $\R$-ext (see \S1.1). This
induces a group isomorphism
$$
\begin{aligned}
\autfun(\A) (\cS \otimes_\R \cS) & \simeq \autfun(\A)(\cS) \times
\cdots \times \,\autfun(\A)(\cS)\\
u &\mapsto (u_\gamma  )_{\gamma  \in\Gamma  }
\end{aligned}
$$
The cocycle condition $u\in Z^1\big(S/R,\,\autfun(\A)\big)$
translates, just as in the classical situation, into the usual
cocycle condition $u_{\gamma  \rho} = u_\gamma  \,^\gamma  u_\rho$
where $\Gamma  $ acts on $\autfun (\A)(\cS) = \,\Aut_{\cS} (\A
\otimes _{\R}\, \cS)$ by conjugation, i.e., $^\gamma  \theta =
(1\otimes \gamma  )\theta    (1\otimes \gamma  ^{-1}).$  This leads
to a natural isomorphism
$$
H^1\big(\cS/\R, \,\autfun (\A)\big) \simeq H^1\big(\Gamma
,\Aut_{\cS}(\A\otimes _{\R}\,\cS)\big)
$$
where the right-hand side is the usual Galois cohomology.}
\end{remark}

\subsection{Limits}\label{limitsss}


{\it Throughout this section, $R:= \K[t^{\pm 1}]$, $S_m:=\K[t^{\pm
1/m}]$, $\widehat{S}:= \us{\longrightarrow}\lim\,S_m, $ and
$\delta_t:=\frac{d}{dt}.$}
\medskip

We are interested in classifying all twisted loop algebras of a
given conformal superalgebra $\A$ over $k.$  If $\sigma  \in
\,\Aut_{k-conf}(\A)$ is of period $m$, then $L(\A,\sigma  )$ is
trivialized by the extension $\cS_m/\R$, where $\R=(R,\delta_t)$ and
$\cS_m=(S_m,\delta_t)$.  (See Example \ref{loopalgdef}.) To compare
$L(\A,\sigma )$ with another $L(\A,\sigma  ^\prime)$ where $\sigma
^\prime$ is of period $m^\prime,$ we may consider a common
refinement $\cS_{mm'}\,.$  An elegant way of taking care of all such
refinements at once is by considering the limit $\widehat{\cS}   =
(\widehat{S},\delta_t).$ For algebras over $k$, and under some
finiteness assumptions, $\widehat{S}  $ plays the role of the
separable closure of $R$ (see \cite{GiPi1} and \cite {GiPi3} for
details). We follow this philosophy in the present situation.

Let $m \in \Z_+$, and let ${\ }^{\overline{\ \ }} \,\, : \Z \to
\Z/m\Z$ be the canonical map. Each extension $\cS_m/\R$ is Galois
with Galois group $\Z/m\Z,$ where $^{\ol 1}(t^{1/m}) = \xi    _m
t^{1/m}$.  (Our choice of roots of unity ensures that the action of
$\Z/\ell m\Z$ on $S_m$ is compatible with that of $\Z/m\Z$ under the
canonical map $\Z/\ell m\Z \to \Z/m\Z.)$  Fix an algebraic closure
$\ol{k(t)}$ of $k(t)$ containing all of the rings $S_m,$ and let
$\pi_1(R)$ be the algebraic fundamental group of $\Spec(R)$ at the
geometric point $a=\,\Spec\,\big(\,\ol{k(t)}\,\big)$.  (See
\cite{SGA1} for details.) Then $\widehat{S}$ is the algebraic
simply-connected cover of $R$, and $\pi_1(R) =
\widehat{\Z}:=\us{\longleftarrow}\lim\,\Z/m\Z,$ where $\widehat{\Z}$
acts continuously on $\widehat{S} $ via $^1 t^{p/q} = \xi ^p_q
t^{p/q}.$  Let $\pi  :\pi_1(R)\to \,\Aut_R(\widehat{S}  )$ be the
corresponding group homomorphism.

Let $\A$ be an $\R = (R,\delta  _t)$-conformal superalgebra. As in
Remark \ref{3.17}, $\pi_1(R)$ acts on $\autfun (\A)(\widehat{\cS}  )
= \;\Aut_{\widehat{\cS} }(\A\otimes_{\R}\widehat{\cS}  )$ by means
of $\pi.$  That is, if $\gamma  \in \pi_1(R)$ and $\theta \in
\,\Aut_{\widehat{\cS}} (\A\otimes_{\R}\, \widehat{\cS}  ),$ then
$$
^\gamma  \theta  = \big(1\otimes \pi (\gamma  )\big) \theta
\big(1\otimes \pi  (\gamma  )^{-1}\big).
$$

Let $\A$ be an $\R$-conformal superalgebra.  Given another
$\R$-conformal superalgebra $\B$, we have
$$
\begin{aligned}
\A \otimes_{\R}\, \cS_m\simeq _{\cS_m}\,\B \otimes _{\R} \cS_{m}
&\Longrightarrow A\otimes_{\R}\cS_{m}\otimes_{\cS_{\bm}} \cS
\simeq_{\cS}
\B \otimes_{\R}\, \cS_m \otimes_{\cS_m} \,\cS\\
&\Longrightarrow A\otimes_{\R} \,\cS \simeq_{\cS} B\otimes_{\R}\,\cS
\end{aligned}
$$
for all extensions $\cS$ of $\cS_m.$  This yields inclusions
$$
H^1\big(\cS_m/\R , \autfun(\A)\big) \subseteq
H_1(\cS_n/\R,\autfun(\A)\big) \subseteq H^1(\widehat{\cS}
/\R,\autfun(\A)\big)
$$
for all $m \vert n,$ hence a natural injective map
\begin{equation}\label{injectivemap}
\eta:\ \us{\longrightarrow}\lim\,H^1(\cS_m/\R, \autfun (\A)\big) \to
H^1(\widehat{\cS}  /\R,\autfun(\A)\big).
\end{equation}

In the classical situation, namely when $\A$ is an algebra, the
surjectivity of the map $\eta$ is a delicate problem (see
\cite{Margaux} for details and references). The following result
addresses this issue for an important class of conformal
superalgebras.

\vskip.3cm
\begin{proposition}\label{limit} Assume that $\A$ satisfies the following
finiteness condition: \vskip.3cm
\begin{description}
\item  {\rm ({\bf Fin})} There exist $a_1,\dots,a_n\in \A$ such that the
set $\{\partial^\ell_{\A}(ra_i)\ |\ r \in R,\,\ell\ge 0\}$ spans
$\A.$
\end{description}

\noindent Then the natural map $\eta \, : \,
\us{\longrightarrow}\lim\,H^1(\cS_m/\R, \autfun (\A)\big) \to
H^1(\widehat{\cS}  /\R,\autfun(\A)\big)$ is bijective. Furthermore,
the profinite group $\pi_1(R)$ acts continuously on
$\Aut_{\widehat{\cS}  } (\A\otimes_{\R}\,\widehat{\cS}  ) =
\,\autfun(\A)(\widehat{\cS}  )$ and
$$
H^1\big(\widehat{\cS}  /\R,\,\autfun(\A)\big) \simeq
H_{\ct}^1\big(\pi_1(R),\,\autfun\,(\A)(\widehat{\cS}  )\big),
$$
where the right $H_{\ct}^1$ denotes the continuous non-abelian
cohomology of the profinite group $\pi_1(R)$ acting (continuously)
on the group $\Aut\,(\A)(\widehat{\cS})$
\end{proposition}

\noindent
\begin{proof} We must show that every $\widehat{\cS}  $-conformal superalgebra isomorphism
$$
\psi  :\A\otimes_{\R} \,\widehat{\cS}  \to \B\otimes_{\R}\,
\widehat{\cS}
$$
is obtained by base change from an $\cS_m$-isomorphism
$$
\psi  _m:\A\otimes_{\R} \,\cS_m \to \B\otimes_{\R}\, \cS_m.
$$
Let $m>0$ be sufficiently large so that $\psi  (a_i\otimes 1) \in
\B\otimes_{\R}\cS_m$ for all $i$.  Since
$\partial_{\B\otimes_{\R}\,\widehat{\cS}  } =
\partial_{\B}\otimes 1+ 1 \otimes \delta_t$ stabilizes
$\B\otimes_{\R}\cS_m$, we have
$$
\begin{aligned}
\psi  \big(\partial^\ell_{\A}(ra_i)\otimes 1\big) &= \psi
\big(\partial^\ell_{\A\otimes_{\R}\, \widehat{\cS}  } (ra_i\otimes
1)\big)\\
&= \partial^\ell_{\B\otimes _{\R}\, \widehat{\cS}  } \big(\psi
(ra_i\otimes 1)\big) = \partial^\ell_{\B\otimes_{\R}\, \widehat{\cS}
} \big(\psi
(a_i\otimes r)\big)\\
&\in \partial^\ell_{\B\otimes_{\R}\, \widehat{S}  } (\B\otimes_\R\,
S_m)\subseteq \B\otimes_R \,\cS_m\,.
\end{aligned}
$$
Thus $\psi  (\A\otimes 1)\subseteq \B\otimes_{\R}\,\cS_m.$  Since
$\psi  $ is $\cS_m$-linear, we get
$$
\psi  (\A\otimes_{\R}\cS_m) \subseteq \B \otimes_R\cS_m\,.
$$
By restriction, we then have an $\cS_m$-conformal superalgebra
homomorphism
$$
\psi  _m: \A\otimes _{\R}\cS_m \to \B\otimes_{\R}\cS_m,
$$
which by base change induces $\psi  .$  Since the extension
$\widehat{S} /S_m$ in $k$-alg is faithfully flat and $\psi  $ is
bijective, our map $\psi  _m$ (viewed as an $S_m$-supermodule map)
is also bijective (by faithfully flat descent for modules). Thus
$\psi  _m$ is a conformal isomorphism of $\cS_m$-algebras.
\end{proof}

An automorphism $\theta  $ of the $\widehat{\cS}$-conformal
superalgebra $\A\otimes_{\R}\,\widehat{\cS}$ is determined by its
restriction to $\A\otimes 1.$  Since $\theta  $ commutes with
$\partial_{\A \otimes_{\R}\widehat{\cS}}$ and since
$\partial^\ell_{\A}(ra_i)\otimes 1 =
\partial^\ell_{\A\otimes _{\R}\,\widehat{\cS}} (ra_i\otimes 1),$ we see that $\theta  $ is
determined by its values on the $a_i\otimes 1.$ Choose $m>0$
sufficiently large so that $\theta  (a_i\otimes 1)\in \A\otimes_R
S_m$ for all $1\le i\le n.$  Then $\{\gamma  \in \widehat{\Z} =
\pi_i(R)\ |\ ^\gamma \theta =\theta  \}$ is of finite index in
$\widehat {\Z},$ hence open (as is well known in the case of the
profinite group $\widehat {\Z}$).\qed

\begin{remark}\label{comparison} {\rm We shall later see that all of the conformal superalgebras that interest us do satisfy the above finiteness condition. Computing \\$H_{\ct}^1\big(\pi_1(R),\,\autfun\,(\A)(\widehat{\cS} )\big)$ is
thus central to the classification of forms. The following two
results are therefore quite useful.
\medskip

(1)  {\it If $\bG$ is a linear algebraic group over $k$ whose
identity connected component is reductive, then the canonical map

$$H_{\ct}^1\big(\pi_1(R),\,\bG(\widehat{\cS}  )\big) \to
H_{\et}^1\big(R,\,\bG \big)$$ is bijective.}

\smallskip

(2)  {\it If $\mathfrak G$ is a reductive group scheme over
$\Spec(R)$ then $H_{\et}^1\big(R,\,\mathfrak G\big) = 1.$}
\medskip

The first result follows from Corollary 2.16.3 of \cite{GiPi3},
while (2) is the main result of \cite{pianzola}.}

\end{remark}

\subsection{The centroid trick}\label{centroidtrick}






Analogous to work done for Lie (super)algebras \cite{ABP2,GrPi}, it
is possible to study the more delicate question of $\K$-conformal
isomorphism, as opposed to the stronger condition of $R$-conformal
isomorphism, using a technique that we will call the {\em centroid
trick}.\footnote{This name was suggested by B.N.~Allison to
emphasize the idea's widespread applicability.}  In this section, we
collect some general facts about centroids and the relationship
between these two types of conformal isomorphism. In the next
section, we will apply the results of this section to some
interesting examples.

Except where otherwise explicitly noted, we assume throughout
\S\ref{centroidtrick} that $\R = (R,\delta_R)$ is an arbitrary
object of $\kalgder.$ Recall that if $R = k$, then $\delta_k = 0$.

For any $\R$-conformal superalgebra $\A$, let $\Ctd_\R(\A)$ be the
set
\begin{equation*}
\{\chi\in\End_{R-smod}(\A)\ |\ \chi(a_{(n)}b)=a_{(n)}\chi(b)\
\hbox{for all}\ a,b\in\A,\ n\in\Z_+\},
\end{equation*}
where $\End_{R-smod}(\A)$ is the set of homogeneous $R$-supermodule
endomorphisms $\A\rightarrow\A$ of degree $\ol{0}$.

Recall that for $r \in R$ we use $r_\A$ to denote the homothety $a
\mapsto ra$. By Axiom (CS3), we have $r_\A \in \Ctd_\R(\A).$ Let
$R_\A = \{r_\A : r \in R\}.$ We have a canonical morphism of
associative $k$- (and $R$-) algebras

\begin{equation}\label{centroidhomothety}
 R \to R_\A \subseteq \Ctd_R(\A).
\end{equation}

Via restriction of scalars, our $\R$-conformal superalgebra $\A$
admits a $k$-conformal structure (where, again, $k$ is viewed as an
object of $\kalgder$ by attaching the zero derivation). This yields
the inclusion
\begin{equation}\label{centroids}
\Ctd_{\R}(\A) \subseteq \Ctd_k(\A).
\end{equation}


\begin{lemma}\label{basiccentroid} Let $\A$ and $\B$ be $\R$-conformal superalgebras such that their restrictions are isomorphic $k$-conformal superalgebras.  That is, suppose there is a $k$-conformal superalgebra isomorphism $\phi : \A \to \B$
(where $\A$ and $\B$ are viewed as $k$-conformal superalgebras by
restriction).  Then the following properties hold.

\vskip.3cm
\begin{description}
\item  {\rm(i)} The map $ \chi  \to \phi\chi  \phi^{-1}$ defines an associative $\K$-algebra isomorphism $\Ctd(\phi):\
\Ctd_k(\A)\rightarrow\Ctd_k(\B)$. Moreover, $\phi$ is an
$\R$-conformal superalgebra isomorphism if and only if
$\Ctd(\phi)(r_\A) = r_\B$ for all $r \in R.$

\item  {\rm (ii)} The map $\delta_{\Ctd_k(\A)}: \chi  \mapsto
[\partial_{_\A},\chi  ] =\partial_{_\A}\chi  -\chi  \partial_{_\A} $
is a derivation of the associative $k$-algebra $\Ctd_k(\A)$.
Furthermore, the  diagram

\begin{equation}\label{pctdcd}
\begin{CD}
\Ctd_k(\A)@> \Ctd(\phi)>>\Ctd_k(\B) \\
@V\delta_{\Ctd_k(\A)}VV @VV\delta_{\Ctd_{k}(\B)}V \\
\Ctd_k(\A)@> \Ctd(\phi)>>\Ctd_k(\B).
\end{CD}
\end{equation}
commutes.

\end{description}
\end{lemma}

\noindent
\begin{proof}
(i) This is a straightforward consequence of the various
definitions.

(ii) For any $\chi\in\Ctd_k(\A)$, $a,b\in\A$ and $n\geq 0$,
\begin{align}
\notag [\p_\A,\chi]\left(a_{(n)}b\right)&=(\p_\A\chi-\chi\p_\A)\left(a_{(n)}b\right)\\
\notag&=\p_\A\left(a_{(n)}\chi (b)\right)-\chi\left(\p_\A (a)_{(n)}b+a_{(n)}\p_\A (b)\right)\\
\notag&=\p_\A (a)_{(n)}\chi (b)+a_{(n)}\p_\A(\chi (b))-\p_\A (a)_{(n)}(\chi b)-a_{(n)}\chi(\p_\A (b))\\
\notag&=a_{(n)}[\p_\A,\chi](b).
\end{align}
The commutativity of Diagram (\ref{pctdcd}) is easy to verify.\qed


\end{proof}


For some of the algebras which interest us the most (e.g.~the
conformal superalgebras in \S\ref{examples}),  the natural ring
homomorphisms $R\rightarrow \Ctd_k(\A)$  are isomorphisms.   This
makes the following result relevant.

\begin{proposition} \label{a1}Let $\A_1$ and $\A_2$ be conformal
superalgebras over $\R = (R,\delta_R).$ Assume that $\Aut_k(\R) =
1,$ i.e., the only $\K$-algebra automorphism of $R$ that commutes
with the derivation $\delta_R$ is the identity.  Also assume that
the canonical maps $R\rightarrow \Ctd_k(\A_i)$ are $\K$-algebra
isomorphisms for $i=1,2$.  Then $\A_1$ and $\A_2$ are isomorphic as
$k$-conformal superalgebras if and only if they are isomorphic as
$\R$-conformal superalgebras.
\end{proposition}

\noindent \proof Clearly, if $\phi:\ \A_1\rightarrow{\A_2}$ is an
isomorphism of $\R$-conformal superalgebras, then it is also an
isomorphism of $k$-conformal superalgebras when $\A_1$ and ${\A_2}$
are viewed as $k$-conformal superalgebras by restriction of scalars.

Now suppose that $\phi:\ \A_1\rightarrow{\A_2}$ is a $k$-conformal
isomorphism, and consider the resulting $k$-algebra isomorphism
$\Ctd(\phi):\ \Ctd_k(\A_1)\rightarrow \Ctd_k({\A_2})$ of Lemma
\ref{basiccentroid}. Under the identification $R_{\A_1} = R =
R_{\A_2},$ the commutativity of Diagram (\ref{pctdcd}), together
with $[\p_{\A_1},r_{\A_1}]=\delta_R(r)_{\A_1}$ and
$[\p_{\A_2},r_{\A_2}]=\delta_R(r)_{\A_2}$, yields that $\Ctd(\phi),$
when viewed as an element of $\Aut_{k-alg}(R),$ commutes with the
action of $\delta_R.$ By hypothesis, $\Ctd(\phi)=\hbox{id}_R$, and
therefore $\phi$ is $R$-linear by Lemma \ref{basiccentroid}(i).
Since $\phi$ also commutes with $\p_{\A_1}$ and preserves
$n$-products, it is an $\R$-conformal isomorphism. \qed

\bigskip

\begin{corollary}\label{3.27}
Let $\A_1$ and $\A_2$ be conformal superalgebras over $\R=(R,\gd_t)$
where $R=k[t,t^{-1}]$ and $\gd_t=\frac{d}{dt}$. Assume that the
canonical maps $R\rightarrow \Ctd_k(\A_i)$ are $k$-algebra
isomorphisms, for $i=1,2$. Then
$$\A_1\cong_{k}\A_2\hbox{\quad if and only if\quad}\A_1\cong_{\R}\A_2.$$
\end{corollary}

\noindent
\begin{proof} The only associative
$\K$-algebra automorphisms of $R$ are given by maps $t\mapsto\ga
t^\eps$, where $\ga$ is a nonzero element of $\K$ and $\eps=\pm 1$.
Thus the only $k$-algebra automorphism commuting with the derivation
$\delta_t$ is the identity map, so the conditions of Proposition
\ref{a1} are satisfied.\qed
\end{proof}

\bigskip

\begin{remark}\label{contrast} {\rm The previous corollary is in sharp contrast to the
situation that arises in the case of twisted loop algebras of
finite-dimensional simple Lie algebras, where $k$-isomorphic forms
need not be $R$-isomorphic.  (See \cite{ABP2} and \cite{pianzola}.)
The rigidity encountered in the conformal case is due to the
presence of the derivation $\delta_t.$}
\end{remark}



\subsection{Central extensions}

In this section we assume that our conformal superalgebras are
Lie-conformal, i.e. they satisfy axioms (CS4) and (CS5). Following
standard practice we  denote the $\lambda$-product $a_{\lambda} b$
by $[a_{\lambda} b]$ (which is then called the $\lambda$-bracket.)

A {\em central extension} of a $k$-conformal superalgebra $\A$ is a
$k$-conformal superalgebra $\widetilde{\A}$ and a conformal
epimorphism $\pi:\ \widetilde{\A}\rightarrow\A$ with kernel
$\ker\pi$ contained in the {\em centre}
$Z(\widetilde{\A})=\{a\in\widetilde{\A}\ |\ \lb{a}{b}=0\ \hbox{for
all}\ b\in\widetilde{\A}\}$ of $\widetilde{\A}$.  Given two central
extensions $(\widetilde{\A},\pi)$ and $(\widetilde{\B},\mu)$ of
$\A$, a {\em morphism} (from $(\widetilde{\A},\pi)$ to
$(\widetilde{\B},\mu)$) in the category of central extensions is a
conformal homomorphism $\phi:\ \widetilde{\A}\rightarrow
\widetilde{\B}$ such that $\mu\circ\phi=\pi$.  A central extension
of $\A$ is {\em universal} if there is a unique morphism from it to
every other central extension of $\A$.

\begin{proposition} Let $(\widetilde{\A_i},\pi_i)$ be central
extensions of $k$-conformal superalgebras $\A_i$ with
$\ker\pi_i=Z(\widetilde{\A_i})$ for $i=1,2$.  Suppose
$\widetilde{\psi}:\ \widetilde{\A_1}\rightarrow\widetilde{\A_2}$ is
a $k$-conformal isomorphism.  Then $\A_1\cong_{k-conf}\A_2$.
\end{proposition}

\proof In the category of $k$-vector spaces, fix sections $\gs_i:\
\A_i\rightarrow\widetilde{\A_i}$ of $\pi_i:\
\widetilde{\A_i}\rightarrow\A_i$.  We verify that
$$\psi=\pi_2\circ\widetilde{\psi}\circ\gs_1:\ \A_1\rightarrow\A_2$$
is a $k$-conformal isomorphism.

Note that
\begin{align*}
\pi_1\lb{\gs_1(x)}{\gs(y)}&=\lb{\pi_1\gs_1(x)}{\pi_1\gs_1(y)}\\
& =\lb{x}{y}\\
& =\pi_1(\gs_1\lb{x}{y}),
\end{align*}
so $\gs_1\lb{x}{y}=\lb{\gs_1(x)}{\gs_1(y)}+w(\gl)$ for some
polynomial $w(\gl)$ in the formal variable $\gl$ with coefficients
in $\ker\pi_1=Z(\widetilde{\A_1})$.

Moreover,
$\lb{\widetilde{\psi}(u)}{\widetilde{\psi}(a)}=\widetilde{\psi}\lb{u}{a}=0$,
for all $u\in Z(\widetilde{\A_1})$ and $a\in\widetilde{\A_1}$, so it
is easy to see that $\wtpsi$ restricts to a bijection between
$Z(\widetilde{\A_1})$ and $Z(\widetilde{\A_2})$.  Therefore,
\begin{align*}
\psi\lb{x}{y}&=\pi_2\circ\wtpsi\circ\gs_1\lb{x}{y}\\
&=\pi_2\circ\wtpsi\big(\lb{\gs_1(x)}{\gs_1(y)}\big)\\
&=\pi_2\circ\wtpsi\lb{\gs_1(x)}{\gs_1(y)}\\
&=\lb{\pi_2\circ\wtpsi\circ\gs_1(x)}{\pi_2\circ\wtpsi\circ\gs_1(y)}\\
&=\lb{\psi(x)}{\psi(y)}.
\end{align*}

To see that $\psi$ commutes with the derivations, note that
$\pi_1\big(\p_{\widetilde{\A_1}}(x)\big)=\p_{\A_1}\big(\pi_1(x)\big)$
for all $x\in\A_1$.  Thus
$\pi_1\Big(\p_{\widetilde{\A_1}}\big(\gs_1(x)\big)\Big)=\p_{\A_1}\big(\pi_1\gs_1(x)\big)=\p_{\A_1}(x)$,
so
$\p_{\widetilde{\A_1}}\big(\gs_1(x)\big)=\gs_1\big(\p_{\A_1}(x)\big)+u$
for some $u\in\ker\pi_1$.  Hence
\begin{align*}
\p_{\A_2}\psi(x)&=\p_{\A_2}\pi_2\circ\wtpsi\circ\gs_1(x)\\
&=\pi_2\circ\wtpsi\big(\p_{\widetilde{\A_1}}\gs_1(x)\big)\\
&=\pi_2\circ\wtpsi\circ\gs_1\big(\p_{\A_1}(x)\big)+\pi_2\circ\wtpsi(u)\\
&=\psi\big(\p_{\A_1}(x)\big)
\end{align*}
since $\wtpsi:\ Z(\widetilde{\A_1})\rightarrow Z(\widetilde{\A_2})$.
Hence $\psi:\ \A_1\rightarrow\A_2$ is a homomorphism of
$k$-conformal superalgebras.

The map $\psi$ is clearly injective: if $x\in\ker\psi$, then
$\wtpsi\circ\psi_1(x)\in\ker\pi_2=Z(\A_2)$, so $\gs_1(x)\in Z(\A_1)$
and $x=\pi_1\psi_1(x)=0$.

To see that $\psi$ is surjective, let $y\in\A_2$.  Then let
$x=\pi_1\circ\wtpsi^{-1}\circ\gs_2(y)$.  For any
$\widetilde{x}\in\widetilde{\A_1}$, we have
$\gs_1\pi_1(\widetilde{x})=\widetilde{x}+v$ for some
$v\in\ker\pi_1$.  Thus
\begin{align*}
\psi(x)&=\pi_2\circ\wtpsi\big(\gs_1\circ\pi_1(\wtpsi^{-1}\circ\gs_2(y))\big)\\
&=\pi_2\circ\wtpsi(\wtpsi^{-1}\circ\gs_2(y)+v)
\end{align*}
for some $v\in\ker\pi_1$.  Then
$$\psi(x)=\pi_2\circ\gs_2(y)+\pi_2\circ\wtpsi(v)=y,$$
and $\psi:\ \A_1\rightarrow\A_2$ is surjective.  Hence $\psi$ is an
isomorphism of $k$-conformal superalgebras.\qed

\begin{corollary} \label{uceiso} Suppose that $(\widetilde{\A_i},\pi_i)$ are
universal central extensions of $k$-conformal superalgebras $\A_i$
with $Z(\A_i)=0$ for $i=1,2$.  Then
$$\widetilde{\A_1}\cong_{k-conf}\widetilde{\A_2}\ \ \ \hbox{if and only
if}\ \ \  \A_1\cong_{k-conf}\A_2.$$
\end{corollary}\qed

\section{Examples and applications}\label{examples}

In this section, we compute the automorphism groups of some
important conformal superalgebras.  It is then easy to explicitly
classify forms of these algebras by computing the relevant
cohomology sets introduced in \S\ref{negy}.

For all of this section, we fix the notation
\begin{eqnarray*}
\R = (R,\delta_R)&:=&\left(\K[t,t^{-1}],\frac{d}{dt}\right)\\
\widehat{\cS}=(\widehat{S},\delta_{\widehat{S}}) &:=&\left(\K[t^q\
|\ q\in {\mathbb{Q}}],\frac{d}{dt}\right),
\end{eqnarray*}
where $\K$ is an algebraically closed field of characteristic zero.

By definition, every $k$-conformal superalgebra is a $\Z/2\Z$-graded
module $\A$ over the polynomial ring $\K[\p]$ where $\p$ acts on
$\A$ via $\p_\A$. For the applications below, we work with
$\K$-conformal superalgebras $\A$ which are free
$\K[\p]$-supermodules. That is, there exists a $\Z/2\Z$-graded
subspace $V=V_{\ol{0}}\oplus V_{\ol{1}}\subseteq\A$ so that
$$\A_{\ol{\iota}}=\K[\p]\ot_\K V_{\ol{\iota}}$$
for $\ol{\iota}=\ol{0},\ol{1}$.

The following result is extremely useful in computing automorphism
groups of conformal superalgebras.

\begin{lemma}\label{free} Let $\A=\K[\p]\ot_\K V$ be a $\K$-conformal superalgebra which is a free $\K[\p]$-supermodule.
Let $\cS = (S, \delta_S)$ be an arbitrary object of $\kalgder.$ Then

(i) Every automorphism of the $\cS$-conformal superalgebra $\A \ot_k
\cS$ is completely determined by its restriction to $ V \simeq (1\ot
V)\ot 1 \subseteq \A \ot_k \cS.$

(ii) Assume $\phi:\ V\ot_\K S\rightarrow V\ot_\K S$ is a bijective
parity-preserving $S$-linear map such that $\phi(\lb{v\ot 1}{w\ot
1})=\lb{\phi(v\ot 1)}{\phi(w\ot 1)}$ for all $v,s\in V$. Then there
is a unique automorphism $\widehat{\phi}\in{\bf Aut}(\A)(\cS)$
extending $\phi$.

\end{lemma}
\proof Let $\{v_i\ |\ i \in I\}$ be a $\K$-basis of $V$ consisting
of homogeneous elements relative to its $\Z/2\Z$-grading, and $\{s_j
\ | \ j \in J \}$ a $k$-basis of $S.$ Since $\p_{\A \ot \cS} = \p_A
\ot 1 + 1 \ot \delta_S$ and since the $v_i$ form a $k[\p]$-basis of
$\A$, the set
$$\{\p_{\A \ot \cS}^\ell(v_i\ot s_j)\ |\ i \in I, j \in J, \ell\geq 0\}$$
is a $\K$-basis of $\A \ot_k S$.  Because any $\cS$-automorphism
must commute with $\p_{\A \ot \cS}$, we have no choice but to define
$\widehat{\phi}$ to be the unique $\K$-linear map on the vector
space $\A \ot_k S$ satisfying $\widehat{\phi}(\p_{\A \ot
\cS}^\ell(v_i\ot s_j))=\p_{\A \ot \cS}^\ell(\phi(v_i\ot s_j)).$ In
particular, we have
\begin{equation}\label{freebasis}
\widehat{\phi}(\p_{\A \ot \cS}^\ell(x))=\p_{\A \ot
\cS}^\ell(\phi(x)).
\end{equation}
for all $x \in V \ot_k S.$ We claim that $\widehat{\phi} \in {\bf
Aut}(\A)(\cS).$ It is immediate from the definition that
$\widehat{\phi}$ is invertible, and that it commutes with the action
of $\p_{\A \ot \cS}.$

For any $v,w\in V$, $r,s\in R$, and $n\in\Z_+$,
$$\begin{aligned}
\phi(v\ot r_{(n)}w\ot s)&= s\phi(v\ot r_{(n)}w\ot 1)\\
&= -s\,p(v,w)\sum_{j=0}^\infty(-1)^{n+j}\p_{\A\ot\cS}^{(j)}\phi(w\ot 1_{(n+j)}v\ot r)\\
&=-s\,p(v,w)\sum_{j=0}^\infty(-1)^{n+j}\p_{\A\ot\cS}^{(j)}r\phi(w\ot 1_{(n+j)}v\ot 1)\\
&=-s\,p(v,w)\sum_{j=0}^\infty(-1)^{n+j}\p_{\A\ot\cS}^{(j)}r\phi(w\ot 1)_{(n+j)}\phi(v\ot 1)\\
&=-s\,p(v,w)\sum_{j=0}^\infty(-1)^{n+j}\p_{\A\ot\cS}^{(j)}\phi(w\ot 1)_{(n+j)}r\phi(v\ot 1)\\
&=s\,(r\phi(v\ot 1))_{(n)}\phi(w\ot 1)\\
&=\phi(v\ot r)_{(n)}\phi(w\ot s),
\end{aligned}$$
by (CS4) and (CS3).  Similar arguments using (CS1) and (CS4) show
that for any homogeneous $x,y\in V\ot S$ and $\ell, m,n\in\Z_+$,
$$\begin{aligned}
\widehat{\phi} \big(\p_{\A \ot \cS}^{(m)} (x)_{(n)} \p_{\A \ot \cS}^{(\ell)}(y) \big)&= \p_{\A \ot \cS}^{(m)} \phi(x)_{(n)} \p_{\A \ot \cS}^{(\ell)} \phi(y)\\
&= \wh \phi \big(\p_{\A \ot \cS}^{(m)}(x)\big)_{(n)} \wh
\phi\big(\p_{\A \ot \cS}^{(\ell)}(y)\big),
\end{aligned}$$
so $\widehat{\phi}$ preserves $n$-products.  Finally, to see that
$\wh \phi$ is $S$-linear, we first observe  that
\begin{equation}\label{trick}
 s_{\A \ot \cS} \circ \partial^{(n)}_{\A \ot \cS} =\os
n{\us{i=0}\sum}(-1)^i\partial^{(i)}_{\A \ot \cS}
 \circ \delta_S(s)^{(n - i)}_{\A \ot  \cS}.
 \end{equation}
This follows from repeated use of Axiom (CS2) applied to the
$\cS$-conformal superalgebra $\A \ot _k \cS$ when taking into
account that $\p_{\A \ot \cS} = \p_\A \ot 1 + 1 \ot \delta_S.$

For all $s\in S$ and $x \in V \ot_k S$, we then have

$$
\begin{aligned}
\wh \phi \big(s_{\A \ot  \cS} & \circ
\partial^{(n)}_{\A\otimes \cS} (x)\big) \\
&= \wh \phi \Big(\sum^n_{i=0} (-1)^i \partial^{(n)}_{\A\otimes \cS}
\circ
\delta_S(s)^{(n - i)}_{\A \ot  \cS} (x)\big) \Big) &{\rm (by \, \ref{trick})}\\
&= \sum^n_{i=0}  (-1)^i\partial^{(n)}_{\A\otimes \cS} \wh \phi
\Big(\delta_S(s)^{(n - i)}_{\A \ot  \cS} (x)\big) \Big) &{\rm (by \, definition \, of}\, \widehat{\phi})  \\
&= \sum^n_{i=0}  (-1)^i\partial^{(n)}_{\A\otimes \cS} \circ
\delta_S(s)^{(n - i)}_{\A \ot  \cS} \phi(x)\big)  &(\cS-{\rm linearity \, on}\, V\otimes_k S)\\
&= s_{\A \ot  \cS} \circ \partial^{(n)}_{\A \ot  \cS} \phi(
x)  &{\rm (by \, \ref{trick})} \\
&= s_{\A \ot  \cS} \circ \wh \phi \big(\partial^{(n)}_{\A\otimes
\cS}(x)\big), &{\rm (by \, definition \, of}\, \widehat{\phi})
\end{aligned}
$$
so $\wh\phi$ commutes with the operator $s_{\A\ot\cS}$, and
$\wh\phi$ is thus $S$-linear.\qed

\subsection{Current  conformal superalgebras}

Let $V=V_{\ol{0}}\oplus V_{\ol{1}}$ be a Lie superalgebra of
arbitrary dimension over the field $\K$.  Let $\Curr V$ be the {\em
current conformal superalgebra} $$\Curr (V):=\K[\p]\ot_\K V$$ with
$n$-products defined by the $\gl$-bracket\footnote{See Remark
\ref{lambdaproduct}.} $\lb{v}{w}=[v,w]$, where $[v,w]$ is the Lie
superbracket for all $v,w\in V$. The derivation $\p_{\Curr(V)}$ is
given by the natural action of $\p$ on the $k$-space $\Curr(V).$

\begin{theorem}\label{currentautos} Let $V$ be a Lie superalgebra over $\K$.  Assume that
for each ideal $W$ of $V$, the centre $Z(W)=\{w\in W\ |\ [w,W]=0\}$
is zero.  Let $\A=\Curr (V).$ Then $\gs(V)\subseteq V$ for all
$\gs\in\Aut_{\K-conf}(\A)$.
\end{theorem}
\proof Let $\pi_V:\A\rightarrow V=\K\ot_\K V$ be the projection of
$\A$ onto the first component of the (vector space) direct sum
$$\A=(\K\ot_\K V)\oplus(\p\K[\p]\ot_\K V).$$
Let $\gs_V=\pi_V\gs:\ V\rightarrow V$, and extend $\gs_V$ to $\A$ by
$\K[\p]$-linearity.

To verify that $\gs_V$ is a $\K$-conformal homomorphism, we expand
both sides of the following equation for all $x,y\in V$:
\begin{equation}\label{exp1}
\gs\lb{x}{y}=\lb{\gs(x)}{\gs(y)}.
\end{equation}
The left-hand side expands as
\begin{align}
\label{lhs} \gs\lb{x}{y}=\gs_V\lb{x}{y}+(\gs-\gs_V)\lb{x}{y}.
\end{align}
The right-hand side is
\begin{align}
\notag \lb{\gs(x)}{\gs(y)}&= \lb{\gs_V(x)}{\gs_V(y)}+\lb{\gs_V(x)}{(\gs-\gs_V)(y)}\\
\label{rhs} &\ \ \ +\
\lb{(\gs-\gs_V)(x)}{\gs_V(y)}+\lb{(\gs-\gs_V)(x)}{(\gs-\gs_V)(y)}.
\end{align}
By (\ref{pr1}) and (\ref{pr2}),
$$\lb{\gs_V(x)}{(\gs-\gs_V)(y)}+\lb{(\gs-\gs_V)(x)}{\gs_V(y)}+\lb{(\gs-\gs_V)(x)}{(\gs-\gs_V)(y)}$$
is contained in the space
$$\K[\gl]\ot\p\K[\p]\ot V+\gl\K[\gl]\ot\K[\p]\ot V,$$
as is $(\gs-\gs_V)\lb{x}{y}$.  Therefore, we can apply $\pi_V$ to
the right-hand sides of (\ref{lhs}) and (\ref{rhs}) and then
evaluate at $\gl=0$ to obtain
\begin{equation}
\gs_V\lb{x}{y}=\lb{\gs_V(x)}{\gs_V(y)}.
\end{equation}
Thus $\gs_V:\ \A\rightarrow\A$ is a homomorphism of $\K$-conformal
superalgebras.

In fact, $\gs_V$ is a $\K$-conformal superalgebra {\em isomorphism}.
By the $\K[\p]$-linearity of $\gs_V$, it is sufficient to verify
that its restriction $\gs_V:\ V\rightarrow V$ is bijective.  But
this is straightforward: if $\gs_V(x)=0$, then $\gs(x)=\p (a)$ for
some $a\in\A$, and $x=\p\gs^{-1}(a)$.  But $x\in V$, so $x=0$, and
$\gs_V$ is injective.  Likewise, if $y\in V$, then write
$\gs^{-1}(y)=z+\p b$ for some $z\in V$ and $b\in\A$.  Then
\begin{eqnarray*}
y&=&\gs_V\gs^{-1}(y)\\
&=&\gs_V(z)+\gs_V(\p b)\\
&=&\gs_V(z)+\pi_V(\p \gs(b))\\
&=&\gs_V(z),
\end{eqnarray*}
so $\gs_V$ is also surjective.  Hence $\gs_V:\ \A\rightarrow\A$ is a
$\K$-conformal automorphism.

Therefore the map
\begin{equation}
\tau=\gs_V^{-1}\gs:\ \A\rightarrow\A
\end{equation}
is a $\K$-conformal automorphism.  Note that $\pi_V\tau(x)=x$ for
all $x\in V$.  Since $\gs=\gs_V\tau$ and $\gs_V(V)\subseteq V$,
Theorem \ref{currentautos} will be proven if we show that
$\tau(V)\subseteq V$.

For nonzero $x\in V$ write
\begin{equation}
\tau(x)=\sum_{i=0}^{M(x)}\p^{(i)}v_{ix}
\end{equation}
for some $v_{ix}\in V$, with $v_{M(x),x}\neq 0$.  Define $v_{i0}$ to
be zero for all $i$ and $M(0)=-1$.  Let
\begin{equation}
W=\Span_\K\{v_{M(x),x}\ |\ x\in V\}.
\end{equation}
We claim that $W\subseteq V$ is an ideal of the Lie algebra $V$.
Indeed, for any $x,y\in V$,
\begin{align}
\notag \tau[x,y]&=\tau\lb{x}{y}\\
\notag &=\lb{\tau(x)}{\tau(y)}\\
\notag &=\left\lbrack\sum_{i=0}^{M(x)}\p^{(i)}v_{ix}  {\,}_{{\,}_{\gl}}\, \sum_{j=0}^{M(y)}\p^{(j)}v_{jy}\right\rbrack\\
\label{tri}
&=\sum_{i=0}^{M(x)}\sum_{j=0}^{M(y)}(-\gl)^{(i)}(\p+\gl)^{(j)}[v_{ix},v_{jy}].
\end{align}
The highest power of $\p$ in (\ref{tri}) is $\p^{(M(y))}$.  Since
the indeterminate $\gl$ does not occur in the expression
$\tau[x,y]$, we see that either
\begin{equation}\label{chetyri}
M([x,y])=M(y)
\end{equation}
or else
\begin{equation}\label{pyat}
[v_{0x},v_{M(y),y}]=0.
\end{equation}
If (\ref{chetyri}) holds, then
$v_{M([x,y]),[x,y]}=[v_{0x},v_{M(y),y}]$, so
\begin{equation}\label{shest}
[x,v_{M(y),y}]\in W
\end{equation}
since $v_{0x}=\pi_V\tau(x)=x$.  If (\ref{pyat}) holds, then
(\ref{shest}) holds trivially since $[x,v_{M(y),y}]=0$.  Therefore
$[V,v_{M(y),y}]\subseteq W$ for all $y\in V$, so $W$ is an ideal of
the Lie superalgebra $V$.

Suppose that $M(x)>0$ for some $x$. Comparing the highest powers of
$\gl$ occuring on both sides of (\ref{tri}), we have
$$0=(-\gl)^{(M(x))}\gl^{(M(y))}[v_{M(x),x},v_{M(y),y}]$$
for all $y\in V$. That is, $[v_{M(x),x},v_{M(y),y}]=0$ for all $y\in
V$, so $v_{M(x),x}$ is in the centre $Z(W)$ of the Lie superalgebra
$W$.  But $Z(W)=0$ by hypothesis, so $v_{M(x),x}=0$, a
contradiction. Hence $M(x)=0$ for all nonzero $x\in V$.  Therefore,
$\tau(x)=v_{0x}=\pi_V\tau(x)=x$ for all $x\in V$, so the
$\K[\p]$-linear map $\tau$ is the identity map on $\A$, and
$\gs=\gs_V\tau=\gs_V:\ V\rightarrow V.$\qed

\begin{corollary}\label{uno}Let $V$ be a simple Lie superalgebra over $\K$.  Then
$$\Aut_{\K-conf}(\Curr(V))=\Aut_{\K-Lie}(V),$$
where $\Aut_{\K-Lie} (V)$ is the group of Lie superalgebra
automorphisms of $V$.
\end{corollary}
\proof Lie automorphisms of $V$ extend uniquely to superconformal
automorphisms of $\Curr(V)$ by $\K[\p]$-linearity.  Conversely,
superconformal automorphisms of $\Curr(V)$ restrict to automorphisms
of the Lie superalgebra $V$ by Theorem \ref{currentautos}. These
correspondences are clearly inverse to one another.\qed

\begin{corollary}\label{zwei}
Let $V$ be a Lie superalgebra over $k$ and let $\cS$ be an extension
in $k-\delta alg$ with the property that for every ideal $W$ of
$V\ot S$, the centre $Z(W)$ is zero.\footnote{For example, any
finite-dimensional simple Lie algebra $V$ over $k$ satisfies this
condition with $\cS=\widehat{\cS}$.}  Then
$$\Aut_{\cS-conf}(\Curr(V) \ot_\K \cS)=\Aut_{S-Lie}(V\ot_\K S).$$

\end{corollary}
\proof Every $\cS$-conformal automorphism $\gs$ of
$\Curr(V)\ot_k\cS$ is also a $k$-conformal automorphism via the
restriction functor. By Theorem \ref{currentautos},
$$\gs(V\ot
S)\subseteq V\ot S.$$ The $\lambda$-bracket in the $\cS$-conformal
superalgebra $\Curr(V)\ot_k \cS$ is
$$[v\ot r_{\lambda}w\ot s]=[v,w]\ot rs$$
for all $v,w\in V$ and $r,s\in S$.  That is,
$\Curr(V)\ot_k\cS=\Curr(V\ot_k S)$ as $\cS$-conformal superalgebras.
Then the argument of Corollary \ref{uno} holds in the
$\cS$-conformal context as well.\qed

\bigskip

\begin{remark}\label{formsofcurrents} {\rm Let $V$ be a finite-dimensional simple Lie superalgebra over $k$.
By Theorem \ref{forms} and Proposition \ref{limit},
 the $\R$-isomorphism classes of
$\widehat{\cS}/\R$-forms of the $\R$-conformal superalgebra
$$\Curr (V)\ot_\K \R=\big(\K[\p]\ot_\K V\big)\ot_\K \R$$
are parametrized by
$$H^1(\widehat{\cS}/\R,{\bf Aut}(\Curr (V))) \simeq H_{\ct}^1\big(\pi_1(R),{\bf Aut}(\Curr (V))
(\widehat{\cS})\big),$$
  By
Corollary \ref{zwei}, ${\bf Aut} (\Curr
(V))(\widehat{\cS})=\Aut_{\wh S-Lie}(V\ot \wh S)$,
  a group that is computed in \cite{GrPi2, kobayashi, PeKa}.  For example, if $V=\mathfrak{sl}_2(\K)$, then
  $\Aut_{\wh S-Lie}(V\ot \wh S) = {\bf PGL}_2(\wh S)$, so
$$H^1(\cS/\R,{\bf Aut}(\Curr (V))=H_{\ct}^1\big(\pi_1(R), {\bf PGL}_2(\wh S)\big)= H_{\et}^1\big(R, {\bf PGL}_2\big) = \{1\}.$$
In particular, all $\widehat{\cS}/\R$-forms of
$\Curr(\mathfrak{sl}_2(k))\ot_k \R$ are trivial, that is, isomorphic
to $\Curr(\mathfrak{sl}_2(k))\ot_k \R$ as an $\R$-conformal
superalgebra.}
\end{remark}

\subsection{Forms of the $N=2$ conformal superalgebra}
Recall that the classical $N=2$ $k$-conformal superalgebra $\A$ is
the free $k[\partial]$-module $\A = k[\partial]\otimes_k\,V$ where
$V=V_{\ol 0}\oplus V_{\ol 1}$,

$$
\begin{aligned}
V_{\ol 0} &= kL \oplus kJ\\
V_{\ol 1} &= kG^+ \oplus kG^-,
\end{aligned}
$$
with $\lambda$-bracket given by\footnote{These conditions say that
$J$ (respectively, $G^\pm)$ is a primary eigenvector of conformal
weight $1$ (resp., $\frac 32)$ with respect to the Virasoro element
$L.$}

\begin{eqnarray}
\lb{L}{L}&=&(\p+2\lambda)L\\
\lb{L}{J}&=&(\p+\lambda)J\\
\lb{L}{G^{\pm}}&=&(\p+\frac{3}{2}\lambda)G^{\pm}\\
\lb{J}{J}&=&0\\
\lb{J}{G^{\pm}}&=&\pm G^{\pm}\\
\lb{G^+}{G^+}&=&\lb{G^-}{G^-}\ =\ 0\\
\lb{G^+}{G^-}&=&L+\frac{1}{2}(\p+2\lambda)J.
\end{eqnarray}


\begin{proposition}\label{4.19}
Let $\widehat{\A} = A\otimes _k \,\widehat{\cS}  .$ Then

\begin{description}
\item  {\rm (1)} For each $s= \alpha  t^q\in \wh{\cS}^\times   ,$
with $\alpha  \in k^\times$ and $q\in \Q,$ there exists a unique
automorphism $\theta  _s \in\,\Aut_{\widehat{\cS}  } (\widehat \A)$
such that
\begin{eqnarray*}
\theta  _s :&& L \mapsto L + q J\otimes t^{-1}\\
&& J \mapsto J\\
&& G^+\mapsto G^+\otimes s\\
&& G^- \mapsto G^- \otimes s^{-1}.
\end{eqnarray*}

\item  {\rm (2)} There exists a unique automorphism $\omega  \in\,
\Aut_{\widehat{\cS}  }(\widehat\A)$ such that

\begin{eqnarray*}
\omega:&& L\mapsto L\\
&&J\mapsto -J\\
&&G^+\mapsto G^-\\
&&G^- \mapsto G^+.
\end{eqnarray*}

\item  {\rm (3)} The map $s\mapsto \theta  _s$ is a group
isomorphism between $\widehat{S}^\times  $ and the subgroup $\langle
\theta  _s\rangle_{s \in \widehat{S}^\times}$ of
$\Aut_{\widehat{\cS}}(\A \ot_k \widehat{\cS})$ generated by the
$\theta _s.$ This isomorphism is compatible with the action of the
algebraic fundamental group $\pi_1(R).$

\item  {\rm (4)} Let $\Z/2\Z$ act on $\widehat{S}^\times  $ by
$^{\overline{1}} t^{p/q} = t^{-p/q}.$  There exists an isomorphism
of $\pi_1(R)$-groups
$$
\psi  :\widehat{S}^\times  \rtimes \Z/2\Z \to\,\Aut_{\widehat{\cS}
} (\A \otimes_k \widehat{\cS}  )
$$
such that
$$
\psi  (s,\overline{\varepsilon}  )\mapsto \theta  _s\omega
^\varepsilon $$ for all $ s \in \widehat{S}^\times$ and $
\varepsilon  =0,1$.
\end{description}
\end{proposition}

\noindent
\begin{proof} The proofs of (1) and (2) are based on Lemma
\ref{free}. One must check that $\theta_s$ and $\omega$ preserve the
$\lambda$-product of any two elements of $V \otimes 1$. This is done
by direct (tedious) calculations.

\bigskip

(3) This is a straightforward consequence of the various
definitions.

\bigskip

(4) The delicate point is to show that the $(\theta_s)_{s \in
\widehat {S}^\times}$ and $\omega$ generate $\Aut_{\widehat{\cS}  }
(\A \otimes_k \widehat{\cS})$. To do this we start with an arbitrary
element $\sigma \in \Aut_{\widehat{\cS}  } (\A \otimes_k
\widehat{\cS})$ and reason as follows:

\medskip

{\it Step 1:  $\sigma (J \otimes 1) = J \otimes s$ for some $s \in
\widehat{S}$.}

This is again a long and tedious computation.  Briefly, write
\begin{equation}\label{a}
\gs(J\ot 1)=\sum_{j=0}^M \pAS^{(j)}(L\ot r_j)+\sum_{k=0}^N
\pAS^{(k)} (J\ot s_k),
\end{equation}
with $r_M\neq 0$.  Note that $\lb{\gs(J\ot 1)}{\gs(J\ot
1)}=\gs\lb{J\ot 1}{J\ot 1}=0.$

Write $\lb{\gs(J\ot 1)}{\gs(J\ot 1)}$ in the form
$$\sum_i \pAS^i(L\ot u_i)+\sum_m\pAS^m(J\ot v_m).$$
Computing with (\ref{a}) shows that $u_{M+1}\neq 0$.  This is
impossible since $\A\ot\wh\cS$ is a free $k[\pAS]$-module and
$\lb{\gs(J\ot 1)}{\gs(J\ot 1)}=0$.  Therefore, $r_M$ cannot be
nonzero.  That is,
$$\gs(J\ot 1)=\sum_{k=0}^N \pAS^{(k)}(J\ot s_k).$$

Comparing the coefficients of the highest powers of $\gl$ occurring
in the relation
$$\gl\gs(J\ot 1)=\gs\lb{J\ot 1}{L\ot 1}=\lb{\gs(J\ot 1)}{\gs(L\ot 1)}$$
shows that $N=0$.  Hence $\gs(J\ot 1)=J\ot s$ for some
$s\in\wh{\cS}$.

should go back to the argument you used in your original 

If we now apply $\sigma$ to $\lb{J \otimes 1}{L \otimes 1} = \lambda
J \otimes 1$, we obtain

\medskip

{\it Step 2:  There exist $c \in k^\times$ and $s \in \widehat{S}$
such that  $\sigma (J \otimes 1) = J \otimes c$ and $\sigma (L
\otimes 1) = L \otimes 1 + J \otimes s.$}

\smallskip

Next we apply $\sigma$ to $\lb{J \otimes 1}{G^\pm \otimes 1} =
\lambda G^\pm \otimes 1.$ This yields

\medskip

{\it Step 3:  There exist $s^+$ and $s^-$ in $\widehat{S}^\times$
such that either

\smallskip

3(a) $\sigma (J \otimes 1) = J \otimes 1$ and $\sigma ( G^\pm
\otimes 1) = G^\pm \otimes s^\pm,$ or

\smallskip

3(b) $\sigma (J \otimes 1) = -J \otimes 1$ and $\sigma ( G^\pm
\otimes 1) = G^\mp \otimes s^\pm.$}
\smallskip

Next we apply $\sigma$ to $\lb{G^+ \otimes 1}{G^- \otimes 1} = L
\otimes 1 + \frac{1}{2} \widehat{\p}(J \otimes 1) + \lambda J
\otimes 1$ to obtain

\medskip

{\it Step 4: If $\sigma$ is as in Case 3(a) above, then $s^+$ and
$s^-$ are inverses of each other. Furthermore, if we write $s^+ = s
= \alpha t^q$ for some $\alpha \in k^\times$ and $q \in \mathbb Q,$
then $\sigma = \theta_s.$}
\smallskip

To finish the proof, we have to consider the case when $\sigma$ is
as in 3(b). Replacing $\sigma$ by $\sigma \omega$ yields an
automorphism that satisfies 3(a), and we can conclude by Step 4.
\qed
\end{proof}




\begin{theorem} \label{4.20} Let $\A$ be the  classical $N=2$ conformal
superalgebra. Up to $k$-conformal isomorphism, there are exactly two
twisted loop algebras of $\A.$  These are $L(\A,\id)$ and
$L(\A,\omega ).$ Furthermore, any $\wh\cS/\R$-form of $\A$ is
isomorphic to one of these two loop algebras.
\end{theorem}

\noindent
\begin{proof} By Theorem \ref{forms}, Proposition \ref{limit},
and Proposition \ref{4.19}, the $\R$-isomorphism classes of
$\widehat{\cS}/\R$-forms of the $\R$-conformal superalgebra
 $\A$
are parametrized by
$$H^1(\wh \cS/\R,{\bf Aut}(\A)) \simeq H_{\ct}^1\big(\pi_1(R),{\bf Aut}(A)(\widehat{\cS})\big)
\simeq H_{\ct}^1\big(\pi_1(R), \widehat{S}^\times  \rtimes
\Z/2\Z)\big).$$ Consider the split exact sequence of $\pi_1(R)
=\wh{\Z}$-groups
\begin{equation}\label{split2}
1\to \widehat{S}^\times   \to \widehat{S}^\times  \rtimes \Z/2\Z \to
\Z/2\Z \to 1.
\end{equation}
Passing to (continuous) cohomology yields
\begin{equation}\label{Exact2}
H^1_{\ct} (\wh{\Z},\widehat{S}^\times  ) \to H^1_{\ct} (\wh{\Z},
\widehat{S}^\times  \rtimes \Z/2\Z) \os \psi  \to
H^1_{\ct}(\wh{\Z},\Z/2\Z).
\end{equation}
The map $\psi$ admits a section (hence is surjective) because the
sequence (\ref{split2}) is split. Since $\wh{\Z}$ acts trivially on
the (abelian) group $\Z/2\Z$, we have
$$
H^1_{\ct}(\wh{\Z},\Z/2\Z) \simeq \,\Hom_{\ct}(\wh{\Z},\Z/2\Z) \simeq
\Z/2\Z.
$$
Note that $\psi$ maps the cohomology classes of $H^1_{\ct} (\wh{\Z},
{\bf Aut}_{\widehat{\cS}}(\A \ot_k \widehat{\cS}))$ corresponding to
the loop algebras $L(\A,\id)$ and $L(\A,\omega)$ to the two distinct
classes of $H^1_{\ct}(\wh{\Z},\Z/2\Z)$.  To prove Theorem
\ref{4.20}, it is thus enough to show that $\psi$ in (\ref{Exact2})
is bijective,  and that $\wh \cS/\R$-forms of $\A$ are $k$-conformal
isomorphic if and only if they are $\R$-conformal isomorphic.

Let ${\bf G_m} = \Spec(k[ z^{\pm 1}])$ denote the multiplicative
group. Recall that $\autfun(\bG_m)\simeq \Z/2\Z$ where the generator
$\ol 1$ of $\Z/2\Z$ acts on $\Spec(k[z^{\pm 1}])$ via $z\mapsto
z^{-1}.$ We now proceed by exploiting the considerations of Remark
\ref{comparison}. Since $H^1_{\et}(R,\Z/2\Z) \simeq \Z/2\Z$, the
bijectivity of $\psi $ translates into the bijectivity of the
analogue map (also denoted by $\psi  )$ at the \'etale level, namely
$$
H^1_{\et}(R,\bG_m) \to H^1_{\et}(R,\bG_m\rtimes \Z/2\Z) \os \psi \to
H^1_{\et}(R,\Z/2\Z)
$$

Since $H_{\et}^1(R,\bG_m) = \,\Pic\,(R)=1$ our map $\psi  $ has
trivial kernel.  The fibre of $\psi  $ over the non--trivial class
of $H^1(R,\Z/2\Z)$ is measured by the cohomology
$H_{\et}^1(R,{\cal{R}}^1_{S_2/R} (\bG_{m}))$ where
${\cal{R}}^1_{S_2/R} (\bG_{m})$ is the twisted form of the
multiplicative $R$-group $\bG_m$ that fits into the exact sequence
$$
1 \longrightarrow {\cal{R}}^1_{S_2/R} (\bG_{m}) \longrightarrow
{\cal{R}}_{S_2/R} (\bG_{m}) \os{\widehat N}\longrightarrow \bG_m
\longrightarrow \bone,
$$
where $\wt N$ comes from the reduced norm $N$ of the quadratic
extension $S_2/R,$ and ${\cal{R}}$ is the Weil restriction. The
functor of points of the $R$-group ${\cal{R}}^1_{S_2/R} (\bG_{m})$
is thus given by
$$
{\cal{R}}^1_{S_2/R} (\bG_{m})(R^\pr) =\{x\in (S_2\otimes_R\,
R^\pr)^\times : N(x) =1\}.
$$
for all $R' \in R-alg.$ Passing to cohomology on this last exact
sequence yields
$$
{\cal{R}}^1_{S_2/R} (\bG_{m})(R) \os N\rightarrow \bG_m (R)
\rightarrow H^1(R,{\cal{R}}^1_{S_2/R} (\bG_{m})) \rightarrow
H^1(R,{\cal{R}}_{S_2/S}(\bG_{m})).
$$
By Shapiro's Lemma,
$$
H^1(R,{\cal{R}}_{S_2/S}(\bG_{m}) ) =\;\Pic\,(S) = 0.
$$
On the other hand, the norm map
$$
\{x\in S^\times_2: N(x) =1\} \os N\longrightarrow R^\times
$$
is surjective. Thus $H^1(R,{\cal{R}}^1_{S_2/R} (\bG_{m})) = 1$ as
desired.\footnote{One can also see that $H^1(R,{\cal{R}}^1_{S_2/R}
(\bG_{m})) = 1$ directly  by applying Remark \ref{comparison} (2) to
the reductive $R$-group ${\cal{R}}^1_{S_2/R} (\bG_{m})$.}

The above cohomological reasoning shows that $L(\A,\id)$ and
$L(\A,\omega  )$ are nonisomorphic as $\R$-conformal superalgebras,
and they represent the $\R$-isomorphism classes of
$\widehat{S}/\R$-forms of $\A.$ To finish the proof, it suffices to
note that the hypotheses of Corollary \ref{3.27} are satisfied, so
(in this case) $\R$-isomorphism is the same as $k$-isomorphism.
This follows easily from the same argument used in the $N=4$ case in
the proof of Theorem \ref{3.63} below. \qed
\end{proof}

\subsection{Forms of the $N=4$ conformal superalgebra}


We now consider the classical $N=4$ conformal superalgebra $\A.$ We
compute the automorphism group ${\bf Aut}(\A)({\wh \cS})$, from
which the existence of infinitely many non-isomorphic twisted loop
algebras will follow easily from the theory developed in
\S\ref{negy}.

We begin by recalling the definition of the $N=4$ conformal
superalgebra $\A$. Let $\A=\A_\even\oplus\A_\odd$ be the $\K$-vector
space with
\begin{equation}\label{N4finite}
\A_\ib=\K[\p]V_\ib\cong V_\ib\ot_\K\K[\p]
\end{equation}
for $\ib=\even,\odd$, and
\begin{eqnarray*}
V_\even&=&\K L\oplus\bigoplus_{s=1}^3\K J^s\\
V_\odd&=&\bigoplus_{a=1}^2\K G^a\oplus\bigoplus_{b=1}^2\K\ol{G}^b.
\end{eqnarray*}
Let $J^s=\frac12\gs^s$, where $\gs^s$ are the Pauli spin matrices
\begin{eqnarray*}
\gs^1&=&\left(\begin{array}{rr}
0 & 1\\
1 &0\\
\end{array}
\right)\\
\gs^2&=&\left(\begin{array}{rr}
0 & -i\\
i &0\\
\end{array}
\right)\\
\gs^3&=&\left(\begin{array}{rr}
1 & 0\\
0 &-1\\
\end{array}
\right).
\end{eqnarray*}
The space $\A$ is a $(\K,\delta)$-conformal superalgebra, with
multiplication given by

\begin{eqnarray}
\lb{L}{L}&=&(\p+2\gl)L \label{LLmult}\\
\lb{L}{J^s}&=&(\p+\gl)J^s\\
\lb{J^m}{J^n}&=&[J^m,J^n]:=J^mJ^n-J^nJ^m\\ \label{LGmult}
\lb{L}{G^a}&=&(\p+\frac32\gl)G^a\\ \label{LGbarmult}
\lb{L}{\ol{G}^a}&=&(\p+\frac32\gl)\ol{G}^a\\ \label{JGmult}
\lb{J^s}{G^a}&=& -\frac12\sum_{b=1}^3\gs_{ab}^sG^b\\
\label{JGbarmult}
\lb{J^s}{\ol{G}^a}&=&\frac12\sum_{b=1}^3\gs_{ba}^s\ol{G}^b\\
\label{GGmult} \lb{G^a}{G^b}&=&\lb{\ol{G}^a}{\ol{G}^b}=0\\
\label{GGbarmult}
\lb{G^a}{\ol{G}^b}&=&2\delta_{ab}L-2(\p+2\gl)\sum_{s=1}^3\gs_{ab}^sJ^s
\end{eqnarray}
for all $m,n,s\in\{1,2,3\}$ and $a,b\in\{1,2\}$, where $\delta_{ab}$
is the Kronecker delta and $\gs_{ab}^s$ is the $(a,b)$-entry of the
matrix $\gs^s$.  The algebra $\A$ is the $N=4$ conformal
superalgebra described in \cite{kac}.

To compute ${\bf Aut}(\A)({\wh \cS})=\Aut_{{\wh \cS}}(\A\ot_\K {\wh
\cS})$, it is enough (by Lemma \ref{free}) to compute the action of
each automorphism in $\Aut_{{\wh \cS}}(\A\ot_\K {\wh \cS})$ on the
subspace $V\ot_\K 1\subseteq\A\ot_\K {\wh \cS}$. Fix
$\gs\in\Aut_{{\wh \cS}}(\A\ot_\K {\wh \cS})$, and choose $d>0$
sufficiently large so that $\gs(V\ot
1)\subseteq\widehat{\A}:=\A\ot\K[t^{1/d},t^{-1/d}]$.  (Such a $d$
exists since $V$ is finite-dimensional.)  Let $z=t^{1/d}$ and
$\widehat{\delta}=\frac{d}{dt}:\ S_d\rightarrow S_d$ with
$S_d=\K[t^{1/d},t^{-1/d}]=\K[z,z^{-1}]$.

Our first step is to show that $\gs$ restricts to an automorphism of
the superconformal subalgebra $\Curr\mathfrak{sl}_2({\wh
S})=\K[\p]\ot\big(\bigoplus_{s=1}^3\K J^s\ot {\wh S}\big).$  Then we
will be able to apply Corollary \ref{zwei}.

Since superconformal automorphisms preserve $\Z/2\Z$-degree,
$$\gs(J^s\ot 1)=\sum_{j=0}^{M_s}\whp^{(j)}(L\ot r_{sj})+u_s$$
for some $r_{sj}\in \wh S$ and $u_s\in\Curr\mathfrak{sl}_2(\wh S)$.
Assume that $r_{sM_s}$ is nonzero if the sum is nonempty (i.e.~if
$M_s\geq 0$).  Suppose that $M_1\geq 0$ and $M_2\geq 0$.  Comparing
the coefficients of $\gl^{M_1+M_2+1}$ on both sides of the equation
\begin{equation}\label{un}
\gs\lb{J^1\ot 1}{J^2\ot 1}=\gs(J^3\ot 1)
\end{equation}
gives
$$(-\gl)^{(M_1)}\gl^{(M_2)}(2\gl)L\ot r_{1M_1}r_{2M_2}=0.$$
That is, $r_{1M_1}=0$ or $r_{2M_2}=0$, a contradiction.  Hence
$M_1<0$ or $M_2<0$.  Then comparing the coefficients of $L$ on both
sides of (\ref{un}) shows that $$0=\sum_{j=0}^{M_3}\whp^{(j)}(L\ot
r_{3j}),$$ so this sum also must be empty and $M_3<0$.  This
argument can be repeated, replacing (\ref{un}) with
$$\gs\lb{J^2\ot 1}{J^3\ot 1}=\gs(J^1\ot 1)$$
and
$$\gs\lb{J^3\ot 1}{J^1\ot 1}=\gs(J^2\ot 1)$$
to show that $M_1<0$ and $M_2<0$, respectively.

Hence $\gs(J^s\ot 1)\in\Curr\mathfrak{sl}_2(\wh S)$ for all $s$, and
$\gs$ restricts to an $\wh \cS$-conformal automorphism of the
subalgebra $\Curr\mathfrak{sl}_2(\wh S)\subseteq \A\ot \wh\cS$.  By
Corollary \ref{zwei}, $\gs$ is the $\K[\p]$-linear extension of an
$\wh S$-Lie algebra automorphism of $\mathfrak{sl}_2(\wh S)$, so by
\cite{kobayashi}, there is some $Y\in {\rm \bf GL}_2(\wh S)$ such that
\begin{equation}
\gs(J^s\ot 1)=YJ^sY^{-1}
\end{equation}
for $s=1,2,3$.  The units of $\wh S$ are the monomials, so $\det
Y=ct^q$ for some $q\in\Q$ and nonzero $c\in\K$.  Since $\K$ is
assumed to be algebraically closed, $c$ has a square root
$\sqrt{c}\in\K$.  Let $\widehat{Y}=\frac{1}{\sqrt{c}}t^{-q/2}Y$.
Then $\widehat{Y}$ has determinant $1$, and conjugation by
$\widehat{Y}$ has the same effect on $J^s$ as conjugation by $Y$, so
we can assume without loss of generality that $Y\in \SL_2({\wh S})$.

Next we consider the image of $L\ot 1$.  Write
$$\gs(L\ot 1)=\sum_{i\in\Z}P_i(\whp)(L\ot z^i)+w$$
for some polynomials $P_i(\whp)$ in the polynomial ring $\K[\whp]$
and $w\in\Curr\mathfrak{sl}_2({\wh S})$.  Note that $\{i\in\Z\ |\
P_i\neq 0\}$ is nonempty (or else $\gs$ would not be surjective).
Set $N=\max\{i\in\Z\ |\ P_i\neq 0\}$ and $M=\min\{i\in\Z\ |\ P_i\neq
0\}$.  If $N>0$, then comparing the coefficients of $L\ot z^{2N}$ on
both sides of
\begin{equation}\label{deux}
\lb{\gs(L\ot 1)}{\gs(L\ot 1)}=\gs\lb{L\ot 1}{L\ot 1}
\end{equation}
gives
$$P_N(-\gl)P_N(\whp+\gl)(\whp+2\gl)(L\ot z^{2N})=0.$$
That is, $P_N=0$, a contradiction.  Hence $N\leq 0$.

If $N<0$, then comparing the coefficients of $L\ot z^{2M-d}$ in
(\ref{deux}) gives
$$P_M(-\gl)P_M(\whp+\gl)(-M/d)(L\ot z^{2M-d})=0,$$
keeping in mind that $\p L\ot z^{NM}=\whp(L\ot
z^M)-L\ot\wh\delta(z^M)$ with $\wh\delta=\frac{d}{dt}$ and
$z=t^{1/d}$.  Hence $P_M=0$, another contradiction.  Thus $N=M=0$.

Let $P(\p):=P_0(\p)$.  Then comparing the coefficients of $L\ot 1$
in (\ref{deux}) gives
$$P(-\gl)P(\whp+\gl)(\whp+2\gl)(L\ot 1)=P(\whp)(\whp+2\gl)(L\ot 1),$$
so $P(\p)$ is a constant (i.e. a member of $\K$) and $P^2=P$.  Since
$\{i\in\Z\ |\ P_i\neq 0\}$ is nonempty, $P\neq 0$, so $P=1$.  Hence
$$\gs(L\ot 1)=L\ot 1+w$$
for some $w\in\Curr\mathfrak{sl}_2({\wh S})$.

Write $w=\sum_{j=0}^{N'}\whp^{(j)}(w_j)$ for some
$w_j\in\mathfrak{sl}_2({\wh S})$, with $w_{N'}\neq 0$.  Suppose
$N'>0$.  For $u\in\mathfrak{sl}_2({\wh S})=\bigoplus_{s=1}^3\K
J^s\ot {\wh S}$,
$$\lb{u}{L\ot 1}=(1\ot\wh\delta)u+\gl u,$$
so
\begin{eqnarray*}
\gs\big((1\ot\wh\delta)u\big)+\gl\gs(u)&=&\lb{\gs(u)}{\gs(L\ot 1)}\\
&=&(1\ot\wh\delta)\gs(u)+\gl\gs(u)+\sum_{j=0}^{N'}(\whp+\gl)^{(j)}[gs(u),w_j],
\end{eqnarray*}
using the fact that $\gs(u)\in\mathfrak{sl}_2({\wh S})$ (Corollary
\ref{zwei}).  Comparing powers of $\whp$ gives $[\gs(u),w_{N'}]=0$.
This holds for all $u\in\mathfrak{sl}_2({\wh S})$, so $w_{N'}$ is in
the centre $Z(\mathfrak{sl}_2({\wh S}))=0$ of the Lie algebra
$\mathfrak{sl}_2({\wh S})$.  Therefore $w_{N'}=0$, a contradiction.
Hence $w\in\mathfrak{sl}_2({\wh S})$.

Hence we have now proven the following lemma:

\begin{lemma}\label{raz}Let $\A$ be the $N=4$ conformal superalgebra defined above.  Then for any $\gs\in\Aut_{\wh\cS-conf}(\A\ot \wh\cS)$, there is some $Y\in \SL_2({\wh S})$ and some $w\in\mathfrak{sl}_2({\wh S})$ so that
\begin{eqnarray*}
\gs(J^s\ot 1)&=&YJ^sY^{-1}\\
\gs(L\ot 1)&=&L\ot 1+w
\end{eqnarray*}
for all $s\in\{1,2,3\}$.\qed
\end{lemma}

Our next task is to find the value of $w$ in Lemma \ref{raz}.  For
all $u\in\mathfrak{sl}_2({\wh S})$, we have
$$(\whp+\gl)\gs^{-1}(u)-\frac{d}{dt}\gs^{-1}(u)=\lb{L\ot 1}{\gs^{-1}(u)},$$
so
\begin{eqnarray*}
(\whp+\gl)u-\gs\big(\frac{d}{dt}\gs^{-1}(u)\big)&=&\lb{\gs(L\ot 1)}{u}\\
&=&\lb{L\ot 1+w}{u}\\
&=&(\whp+\gl)u-\frac{d}{dt}u+[w,u],
\end{eqnarray*}
and
\begin{eqnarray*}
[w,u]&=&\frac{d}{dt}u-\gs\big(\frac{d}{dt}\gs^{-1}(u)\big)\\
&=&u'-Y(Y^{-1}uY)'Y^{-1},
\end{eqnarray*}
where prime ($'$) denotes the derivative taken with respect to the
variable $t$.  But
\begin{equation}\label{Yeqn}
Y'Y^{-1}+Y(Y^{-1})'=(YY^{-1})'=0,
\end{equation}
 so
\begin{eqnarray*}
[w,u]&=&u'-Y(Y^{-1}uY)'Y^{-1}\\
&=&-Y(Y^{-1})'u-uY'Y^{-1}\\
&=&[Y'Y^{-1},u].
\end{eqnarray*}
Thus $w-Y'Y^{-1}$ is in the centralizer of $\mathfrak{sl}_2({\wh
S})$ in $\mathfrak{gl}_2({\wh S})$.  But writing
$$Y=\left(\begin{array}{cc}
e & f\\
g & h
\end{array}
\right),$$ with $e,f,g,h\in {\wh S}$, a quick computation shows that
the trace $tr(Y'Y^{-1})=(eh-fg)'$.  But $(eh-fg)'=0$ since $Y\in
\SL_2({\wh S})$.  Hence, $w-Y'Y^{-1}\in Z(\mathfrak{sl}_2({\wh
S}))=0$, the centre of $\mathfrak{sl}_2({\wh S})$, so we have the
following proposition.
\begin{proposition}\label{pair}
Let $\A$ be the $N=4$ conformal superalgebra defined above.  Then
for any $\gs\in\Aut_{{\wh \cS}-conf}(\A\ot \wh\cS)$, there is some
$Y\in \SL_2({\wh S})$ so that
\begin{eqnarray*}
\gs(J^s\ot 1)&=&YJ^sY^{-1}\\
\gs(L\ot 1)&=&L\ot 1+Y'Y^{-1}
\end{eqnarray*}
for all $s\in\{1,2,3\}$.\qed
\end{proposition}

Next we consider the action of $\gs$ on the odd part of $\A\ot
\wh\cS$.  Let $v\in V_\odd$.  Write
$$\gs(v\ot 1)=\sum_{i=0}^{M'}\sum_{s=1}^4\whp^{(i)}(v_s\ot r_{si}),$$
where $v_1=G^1$, $v_2=G^2$, $v_3=\ol{G}^1$, $v_4=\ol{G}^2$, and
$r_{sM'}\neq 0$ for some $s$.  Then writing $w$ for $Y'Y^{-1}$, we
have
\begin{eqnarray*}
(\whp+\frac32\gl)\gs(v\ot 1)&=&\gs\lb{L\ot 1}{v\ot 1}\\
&=&\lb{\gs(L\ot 1)}{\gs(v\ot 1)}\\
&=&\sum_{i=0}^{M'}\sum_{s=1}^4(\whp+\gl)^{(i)}\big((\whp+\frac32\gl)(v_s\ot r_{si})\\
&&\ \ \ \ \ \ \ \ \ \ \ \ \ \ \ \ \ \ \ \ \ \ \ \ \ \ -\
v_s\ot\wh\delta(r_{si})+\lb{w}{v_s\ot r_{si}}\big).
\end{eqnarray*}
From the definition of the relevant products, (\ref{LGmult}),
(\ref{LGbarmult}), and (\ref{extproddef}), it is clear that
$\lb{w}{v_s\ot r_{si}}$ contains no nonzero powers of $\gl$.  If
$M'>0$, then comparing the coefficients of $\lambda^{M'+1}$ gives
$$0=\sum_{s=1}^4\frac32(M'+1)\gl^{(M'+1)}v_s\ot r_{sM'}.$$
Thus $r_{sM'}=0$ for all $s$, a contradiction.  Hence $M'\leq 0$,
and $\gs(v\ot 1)\in V\ot S$.  Therefore, by the ${\wh S}$-linearity
of $\gs$ and Proposition \ref{pair}, we have proven the following
lemma.

\begin{lemma} Let $\gs\in\Aut_{\wh\cS}(\A\ot \wh\cS)$.  Then $\gs(V\ot \wh S)\subseteq V\ot \wh S.$\qed
\end{lemma}

For the computations that follow, it will be helpful to use the
following notation:
\begin{eqnarray*}
\left(\begin{array}{c}
a\\
b\end{array} \right)\ot \left(\begin{array}{c}
1\\
0\end{array}
\right):=G^1\ot a+G^2\ot b\\
\left(\begin{array}{c}
a\\
b\end{array} \right)\ot \left(\begin{array}{c}
0\\
1\end{array} \right):=\ol{G}^1\ot a+\ol{G}^2\ot b
\end{eqnarray*}
for all $a,b\in {\wh S}$.  In this notation, the relations
(\ref{JGmult}), (\ref{JGbarmult}), and (\ref{GGbarmult}) become
\begin{eqnarray}
\label{JGmult'} \lb{J^s}{\left(\begin{array}{c}
a\\
b\end{array}\right)\ot \left(\begin{array}{c}
1\\
0\end{array}\right)}&=&-(J^s)^T\left(\begin{array}{c}
a\\
b\end{array}\right)\ot \left(\begin{array}{c}
1\\
0\end{array}\right)\\
\label{JGbarmult'}\lb{J^s}{\left(\begin{array}{c}
g\\
h\end{array}\right)\ot \left(\begin{array}{c}
0\\
1\end{array}\right)}&=&J^s\left(\begin{array}{c}
g\\
h\end{array}\right)\ot \left(\begin{array}{c}
0\\
1\end{array}\right)
\end{eqnarray}
\begin{align}
\notag \lb{\left(\begin{array}{c}
a\\
b\end{array}\right)\ot \left(\begin{array}{c}
1\\
0\end{array}\right)}{\left(\begin{array}{c}
g\\
h\end{array}\right)\ot \left(\begin{array}{c}
0\\
1\end{array}\right)}&\\
=2\left(\begin{array}{cc} a & b \end{array}\right)
\left(\begin{array}{c}
g\\
h\end{array}\right)L &-\
2(\p+2\gl)\sum_{s=1}^3\left(\begin{array}{cc} a &
b\end{array}\right)\gs^s\left(\begin{array}{c}
g\\
h
\end{array}\right)J^s.\label{GGbarmult'}
\end{align}
for all $a,b,g,h\in\K$ and $s\in\{1,2,3\}$, where $(J^s)^T$ is the
transpose of the matrix $J^s$.

Using the fact that $\gs$ preserves the relation (\ref{JGmult'}), we
see that for any $a,b\in\K$,
\begin{align*}
-(YJ^sY^{-1})^T&\gs_1\left(\abG\right)+YJ^sY^{-1}\gs_2\left(\abG\right)\\
=-\gs_1&\left((J^s)^T\abG\right)-\gs_2\left((J^s)^T\abG\right),
\end{align*}
where $\gs_i:=\pi_i\gs$ and $\pi_1$ (resp., $\pi_2$) is the
projection of $\A_\odd\ot {\wh S}$ onto
$$W_1:=\bigoplus_{j=1}^2\K G^j\ot {\wh S}\ \ \ \left(\hbox{resp.,\ }W_2:=\bigoplus_{j=1}^2\K \ol{G}^j\ot {\wh S}\right)$$
in the direct sum
$$\A_\odd=W_1\oplus W_2.$$
Then
\begin{equation}\label{csillag}
-(YJ^sY^{-1})^T\gs_1=-\gs_1(J^s)^T
\end{equation}
as ${\wh S}$-linear maps on the vector space $W_1$.  Let
$v=\left(\begin{array}{r}
v_1\\
v_2
\end{array}\right)\ot\left(\begin{array}{r}
1\\
0
\end{array}\right)$ be in the kernel $\ker\gs_1$ of the restriction $\gs_1:\ W_1\rightarrow W_1$.  Then by (\ref{csillag}),
$$Mv:=\left(M\left(\begin{array}{r}
v_1\\
v_2
\end{array}\right)\right)\ot\left(\begin{array}{r}
1\\
0
\end{array}\right)\in\ker\gs_1$$
for all $M\in{\mathfrak{sl}}_2(\wh S)$.  Thus $\ker\gs_1=0$ or
$\ker\gs_1=W_1$.  If $\ker\gs_1=0$, we see by (\ref{csillag}) that
conjugation by $\gs_1$ has the same effect on $\mathfrak{sl}_2(\K)$
as conjugation by $(Y^{-1})^T$. The centralizer of
$\mathfrak{sl}_2(\K)$ in ${\rm \bf GL}_2({\wh S})$ consists of matrices of
the form $cz^mI$, where $c\in\K^{\times}$, $m\in\Z$, and $I$ is the
$2\times 2$ identity matrix.  Thus
\begin{equation}\label{ketcsillag}
\gs_1=cz^m(Y^{-1})^T:\ W_1\rightarrow W_1
\end{equation}
for some $c\in \K^{\times}$ and $m\in\Z$.  If $\ker\gs_1=W_1$, then
obviously (\ref{ketcsillag}) also holds, with $c=0$.

By a similar argument,
$$\gs_2=dz^nY\left(\begin{array}{rr}
0 & 1\\
-1 & 0
\end{array}
\right):\ W_1\rightarrow W_2$$ for some $d\in\K$ and $n\in\Z$.  Thus
\begin{align}
\notag \gs\left(\abG\right)&\\
\label{abGim}=cz^m(Y^{-1})^T\left(\begin{array}{c}
a\\
b\end{array}\right)&\ot \left(\begin{array}{c}
1\\
0\end{array}\right)+dz^nY\left(\begin{array}{rr}
0 & 1\\
-1 & 0
\end{array}
\right)\left(\begin{array}{c}
a\\
b\end{array}\right)&\ot \left(\begin{array}{c}
0\\
1\end{array}\right).
\end{align}

Repeating this argument on $W_2$ using relation (\ref{JGbarmult'}),
we see that for some fixed $e,f\in\K$ and $k,\ell\in\Q$,
\begin{align}
\notag \gs\left(\ghGbar\right)&\\
=ez^k(Y^{-1})^T&\left(\begin{array}{rr}
0 &1\\
-1 &0
\end{array}
\right)\left(\begin{array}{c}
g\\
h\end{array}\right)&\ot \left(\begin{array}{c}
1\\
0\end{array}\right)+fz^\ell Y\ghGbar
\end{align}
for all $g,h\in\K$.

From (\ref{LGmult}), we see that
\begin{equation}\label{ket}
\lb{L\ot 1}{\abG}=(\whp+\frac32\gl)\abG
\end{equation}
for all $a,b\in\K$.  Apply $\gs_1$ to both sides of (\ref{ket}) and
compute using (\ref{abGim}).  Equating the terms on both sides of
the equation which are constant with respect to $\gl$ and $\whp$
gives
$$(cz^m(Y^{-1})^T)'=-(Y'Y^{-1})^Tcz^m(Y^{-1})^T,$$
where prime ($'$) denotes element-by-element differentiation with
respect to $t=z^d$. Taking the transpose of both sides and
simplifying using (\ref{Yeqn}) gives $c=0$ or $mz^{m-1}Y^{-1}=0$.
That is, $c=0$ or $m=0$.  If $c=0$, then we can obviously assume
that $m=0$.  Similarly, $n=k=\ell=0$.

By (\ref{GGbarmult'}), the following equation holds for all
$a,b,g,h\in\K$:
\begin{align}
\notag \left\lbrack\gs\left(\abG\right)_\gl\right.&\left.\gs\left(\ghGbar\right)\right\rbrack\\
\label{anotherannoyingeqn}=2\left(\begin{array}{cc} a & b
\end{array}
\right)\left(\begin{array}{r}
g\\
h
\end{array}\right)\gs&(L\ot 1)-2(\whp+2\gl)\gs\left(\sum_{s=1}^3\left(\begin{array}{cc}
a & b
\end{array}
\right)\gs^s\left(\begin{array}{r}
g\\
h
\end{array}\right)\right)J^s.
\end{align}

  Then comparing the coefficients of $L\ot 1$ on both sides of (\ref{anotherannoyingeqn}) gives
$$2cf\left(\begin{array}{cc}
a & b
\end{array}
\right)\left(\begin{array}{r}
g\\
h
\end{array}\right)-2de\left(\begin{array}{cc}
g & h
\end{array}
\right)\left(\begin{array}{r}
a\\
b
\end{array}\right)=2\left(\begin{array}{cc}
a & b
\end{array}
\right)\left(\begin{array}{r}
g\\
h
\end{array}\right),$$
so $cf-de=1$.  Thus we have the following proposition:
\begin{proposition}\label{autosprop}
Let $\gs\in\Aut_{\wh\cS}(\A\ot \wh\cS)$. Then for some $Y\in
\SL_2({\wh S})$ and $\left(\begin{array}{rr}
c & d\\
e & f
\end{array}
\right)\in \SL_2(\K)$, $\gs$ satisfies the following formulas
\begin{align}
\label{ODIN} \gs(L\ot 1)&=L\ot 1+Y'Y^{-1}\\
\label{DVA}\gs(J^s\ot 1)&=YJ^sY^{-1}\\
\notag \gs\left(\abG\right)&=c(Y^{-1})^T\abG\\
\label{TRI}&\quad\quad +dY\left(\begin{array}{rr}
0 & 1\\
-1 & 0
\end{array}\right)\left(\begin{array}{c}
a\\
b\end{array}\right)\ot\left(\begin{array}{c}
0\\
1\end{array}\right)\\
\notag \gs\left(\left(\begin{array}{c}
a\\
b\end{array}\right)\ot\left(\begin{array}{c}
0\\
1\end{array}\right)\right)&=e(Y^{-1})^T\left(\begin{array}{rr}
0 & 1\\
-1 & 0
\end{array}\right)\abG\\
\label{chetyre}&\quad\quad + fY\left(\begin{array}{c}
a\\
b\end{array}\right)\ot\left(\begin{array}{c}
0\\
1\end{array}\right)
\end{align}
for all $a,b\in\K$ and $s=1,2,3$.
\end{proposition}\qed

\bigskip

The converse to Proposition \ref{autosprop} is that given any $Y\in
\SL_2({\wh S})$ and $\left(\begin{array}{rr}
c & d\\
e & f\end{array}\right)\in \SL_2(\K)$, (\ref{ODIN})--(\ref{chetyre})
defines an automorphism $\gs\in\Aut_{\wh\cS}(\A\ot \wh\cS)$.  This
follows (using Lemma \ref{free}) from the long and tedious
verification that the ${\wh S}$-linear map $\gs:\A\ot
\wh\cS\rightarrow\A\ot \wh\cS$ defined by
(\ref{ODIN})--(\ref{chetyre}) preserves the $\lambda$-bracket on the
following relation:
\begin{equation}\label{PYAT}
\gs\lb{w_1\ot 1}{w_2\ot 1}=\lb{\gs(w_1\ot 1)}{\gs(w_2\ot 1)}
\end{equation}
for $w_1,w_2\in\{L,J^s,G^i,\ol{G}^i\ |\ s=1,2,3,\ i=1,2\}$.

To determine the group structure on the set of automorphisms, fix
$Y\in \SL_2({\wh S})$ and $X\in \SL_2(\K)$.  Let $\eta_1$ (resp.,
$\eta_2$) denote the automorphism determined by $(Y,I)$ (resp.,
$(I,X)$), where $I$ is the $2\times 2$ identity matrix.  Then it is
straightforward to verify that
$$\eta_2\eta_1\eta_2^{-1}=\eta_1,$$
so there is a group epimorphism
$$\phi: \SL_2({\wh S})\times \SL_2(\K)\rightarrow \Aut_{\wh\cS}(\A\ot \wh\cS).$$

Finally, suppose that the automorphism determined by the pair
$(Y,X)\in \SL_2({\wh S})\times \SL_2(\K)$ is in the kernel
$\ker\phi$. Writing $X=\left(\begin{array}{cc}
c & d\\
e & f\end{array}\right)$,
$$c(Y^{-1})^T\abG=\abG$$
for all $a,b\in\K$ by (\ref{TRI}).  Thus $c(Y^{-1})^T=I$, so $Y=cI$
and $c=\pm 1$.  By (\ref{TRI}) and (\ref{chetyre}), we also see that
$d=e=0$ and $Y=fI$.  Since $(-I,-I)$ determines the identity map on
$\A\ot \wh\cS$ by (\ref{ODIN})--(\ref{chetyre}), the kernel of
$\phi$ is the  subgroup of $\SL_2({\wh S})\times \SL_2(\K)$
generated by $(-I,-I)$:
$$\ker(\phi) =\langle(-I,-I)\rangle \simeq \Z/2\Z$$
We have now proven the following:
\begin{proposition}\label{auto4}
Let $\A$ be the $N=4$ conformal superalgebra defined above, and let
${\wh S}=\K[t^q\ |\ q\in\Q]$.  Then
\begin{equation}
\Aut_{\wh\cS-conf}(\A\ot_\K \wh\cS)=\frac{\SL_2(\wh S)\times
\SL_2(\K)}{\langle(-I,-I)\rangle}.\end{equation}
\end{proposition} \qed

Applying our theory of forms ( \S 2) now allows us to classify
twisted loop algebras of the $N=4$ conformal superalgebra $\A$.

\begin{theorem}\label{3.63} Let $\A$ be the $N=4$ conformal superalgebra defined above.
Then there are canonical bijections between the following sets:

\begin{description}
\item  {\rm (i)} $\R$-isomorphism classes of twisted loop algebras of $\A$,

\item  {\rm (ii)} $k$-isomorphism classes of twisted loop algebras of $\A$,

\item  {\rm (iii)} $\R$-isomorphism classes of $\wh\cS/\R$-forms of
$\A\ot_k\R$,

\item  {\rm (iv)} conjugacy classes of elements of finite order
in ${\bf PGL}_2(\K)$,

\end{description}
\end{theorem}

\noindent
\begin{proof}  Let $\A_m$ be the $\cS_m$-conformal superalgebra $\A\ot_k\cS_m$
where $\cS_m=(S_m,\frac{d}{dt})$, $S_m=k[t^{\pm\frac{1}{m}}]$, and
$m\geq 1$. Let $V=V_\even\oplus V_\odd$, where $V_\even$ and
$V_\odd$ are defined as in (\ref{N4finite}).  We divide our proof
that (i) and (ii) are equivalent into several steps.

\bigskip

{\em Step 1:\ $\A_m=\Span_k\{v_{(1)}\p_\A^{(\ell)}L\ot 1\ |\ v\in
V\ot S_m,\ \ell\geq 0\}$}

{\em Proof:}  We show that $\p_\A^{\ell}V\ot
S_m\subseteq\sum_{j=0}^\ell V\ot {S_m}_{(1)}\p_\A^{(j)}L\ot 1$ using
induction on $\ell$.  For $\ell=0$, we see that $V\ot
{S_m}_{(1)}L\ot 1=\big(V_{(1)}L\big)\ot S_m=V\ot S_m$, since
$L,J^s,G^i,\overline{G}^i$ are {\em primary eigenvectors} of $L$
with conformal weight greater than $1$ for $s=1,2,3$ and $i=1,2$.
That is,
\begin{align*}
a_{(0)}L&=(\Delta-1)\p_\A a\\
a_{(1)}L&=\Delta a\\
a_{(m)}L&=0
\end{align*}
for all $m>1$, $a=L,J^s,G^i,\overline{G}^i$, and some
$\Delta=\Delta(a)\geq 1$.

It is straightforward to verify that for $\ell\geq 1$ and $s\in
S_m$,
$$ a\ot s_{(1)}\p_\A^{(\ell+1)}L\ot 1=\big(
(\ell+2)\Delta-(\ell+1)\big)\p_\A^{(\ell+1)}a\ot
s+\Delta\p_\A^{(\ell)}a\ot\frac{ds}{dt}.$$ Since $\Delta\geq 1$, we
see that $\Delta\p_\A^{(\ell+1)}a\ot s\in\sum_{j=0}^{\ell+1}V\ot
S_m,\ \ell\geq 0\}$.

\bigskip

{\em Step 2:\ Let $\B=\Loop(\A,\gs)\subseteq\widehat{\A}$ for some
finite order automorphism $\gs:\ \A\rightarrow\A$.  Then
$\B=\Span_k\{a_{(1)}\p_\A^{(\ell)}L\ot 1\ |\ a\in\B,\ \ell\geq
0\}$.}

{\em Proof:}  Let $\Gamma\subseteq\Aut_{\cS_m}\A_m$ be the cyclic
subgroup of order $m:=|\gs|$ generated by $\gs\ot\psi$, where $\psi$
is the $\R$-automorphism of $\cS_m$ given by sending $t^{\frac1m}$
to $\xi_m^{-1}t^{\frac{1}{m}}$ and $\xi_m$ is the primitive $m$th
root of $1$ fixed in \S0.  Let
\begin{eqnarray*}
\pi:\ \A_m&\rightarrow&\A_m\\
a&\mapsto&\frac{1}{m}\sum_{i=0}^{m-1}(\gs\ot\psi)^i(a).
\end{eqnarray*}
Then $\B=\A_m^\Gamma$, the set of $\Gamma$-fixed points in $\A_m$,
and $\pi$ is a surjection from $A_m$ to $\B$.

Since $\A_m=\Span_k\{v_{(1)}\p_\A^{(\ell)}L\ot 1\ |\ v\in V\ot S_m,\
\ell\geq 0\}$ and $\gs\ot \psi\in\Aut_{\cS_m}(\A_m)$ by Lemma
\ref{preservesnprod}, we have
\begin{align*}
\B=\pi(\A_m)&=\Span_k\left\{\pi\left(v_{(1)}\p_\A^{(\ell)}L\ot 1\right)\ |\ v\in V\ot S_m,\ \ell\geq 0\right\}\\
&=\Span_k\left\{\sum_{i=0}^{m-1}\pi(v)_{(1)}\p_\A^{(\ell)}\gs^i(L)\ot
1\ |\ v\in V\ot S_m, \ \ell\geq 0\right\}.
\end{align*}
By \ref{TRI}, there is a $Y\in \SL_2(S_m)$ such that
\begin{align*}
\gs(G^1)\ot 1&=(\gs\ot 1)(G^1\ot 1)\\
&=c(Y^{-1})^T\abG +dY\left(\begin{array}{rr}
0 & 1\\
-1 & 0
\end{array}\right)\left(\begin{array}{c}
1\\
0\end{array}\right)\ot\left(\begin{array}{c}
0\\
1\end{array}\right),
\end{align*}
so $Y\in \SL_2(k)$.  Then $Y'=0$, so $\gs(L)\ot 1=(\gs\ot 1)(L\ot
1)=L\ot 1$ by (\ref{ODIN}).  Hence $\gs(L)=L$.

Therefore,
\begin{align*}
\B&=\Span_k\left\{\sum_{i=0}^{m-1}\pi(v)_{(1)}\p_\A^{(\ell)}\gs^i(L)\ot 1\ |\ v\in V\ot S_m, \ \ell\geq 0\right\}\\
&=\Span_k\left\{\pi(v)_{(1)}\p_\A^{(\ell)}L\ot 1\ |\ v\in V\ot S_m, \ \ell\geq 0\right\}\\
&\subseteq\Span_k\{a_{(1)}\p_\A^{(\ell)}L\ot 1\ |\ a\in\B,\ \ell\geq
0\}.\
\end{align*}

\bigskip

{\em Step 3:\ Let $\chi\in\Ctd_k(\B)$, where $\B$ is as above.  Then
$\chi(L\ot 1)=L\ot r$ for some $r\in R$.}

{\em Proof:}  By the argument in Step 2, every automorphism
$\gs\in\Aut_{k-conf}\A$ fixes $L$, so $L\ot 1\in\B$.  Then
$$L\ot 1_{(1)}\chi(L\ot 1)=\chi(L\ot 1_{(1)}L\ot 1)=2\chi(L\ot 1).$$
Taking an eigenspace decomposition of $\widehat{\A}$ with respect to
the operator
$$L\ot 1_{(1)}:\ a\ot s\mapsto L\ot 1_{(1)}a\ot s=\left(L_{(1)}a\right)\ot s,$$
we have
\begin{equation}\label{espacedecomp}
\widehat{\A}=\left(\bigoplus_{k=1}^\infty\widehat{\A}_k\right)\oplus\left(\bigoplus_{\ell=1}^\infty\widehat{\A}_{\frac{1}{2}+\ell}\right),
\end{equation}
where
\begin{eqnarray*}
\widehat{\A}_k&=&\Span\left\{\p_\A^{(k-2)}L\ot r,\p_\A^{(k-1)}J\ot r\ \big|\ J\in\{J^1,J^2,J^3\}\ \hbox{and}\ r\in\widehat{S}\right\}\\
\widehat{\A}_{\frac{1}{2}+\ell}&=&\Span\left\{\p_\A^{(\ell-1)}G\ot
r\ \big|\ G\in\{G^1,G^2,\overline{G}^1,\overline{G}^2\}\ \hbox{and}\
r\in\widehat{S}\right\}
\end{eqnarray*}
are the eigenspaces with eigenvalues $k$ and $\frac{1}{2}+\ell$,
respectively.

Thus $\chi(L\ot 1)\in\widehat{\A}_2$, so $\chi(L\ot 1)=L\ot
r+\sum_{i=1}^3\p_\A J^i\ot r_i$ for some $r,r_i\in\widehat{S}$.  But
also
\begin{align*}
0&=\chi(L\ot 1_{(2)}L\ot 1)\\
&=L\ot 1_{(2)}\chi(L\ot 1)\\
&=L\ot 1_{(2)}\left(L\ot r+\sum_{i=1}^3\p_\A J^i\ot r_i\right)\\
&=2\sum_{i=1}^3J^i\ot r_i.
\end{align*}
Hence $r_i=0$ for $i=1,2,3$, so $\chi(L\ot 1)=L\ot r$.  Since
$\chi\in\Ctd_k(\B)$ and $\gs(L)=L$, we see that $r\in R$.

\bigskip

{\em Step 4:\ Let $\chi$ and $r$ be as in Step 3.  Then
$\chi\left(\p_\A^{(k)}L\ot 1\right)=\p_\A^{(k)}L\ot r$ for all
$k\geq 0$.}

{\em Proof:}  If $b\in\B$ is a primary eigenvector of itself with
conformal weight $\Delta$, then it is straightforward to verify that
$b_{(k+1)}\p_\B^{(k)}b=(\Delta+k-1)\p_\B^{(k)}b$.
By (CS4), we have
\begin{align*}
(\Delta+k-1)\chi\left(\p_\B^{(k)}b\right)&=\chi\left(b_{(k+1)}\p_\B^{(k)}b\right)\\
&=(-1)^k\chi\left(\p_\B^{(k)}b_{(k+1)}b\right)\\
&=(-1)^k\p_\B^{(k)}b_{(k+1)}\chi(b).
\end{align*}
If $\chi(b)=r b$, this gives
\begin{align*}
(\Delta+k-1)\chi\left(\p_\B^{(k)}b\right)&=r_\B(-1)^k\p_\B^{(k)}b_{(k+1)}b\\
&=r_\B b_{(k+1)}\p_\B^{(k)}b\\
&=(\Delta+k-1)r\p_\B^{(k)}b,
\end{align*}
so $\chi\left(\p_\B^{(k)}b\right)=r_{\B} \circ \p_\B^{(k)}b$ if
$\Delta>1$. Applying this result to $b=L\ot 1$, we see that
\begin{align*}
\chi\left(\p_\A^{(k)}L\ot 1\right)&=\chi\left(\p_\B^{(k)}(L\ot 1)\right)\\
&=r_\B \circ \p_\B^{(k)}(L\ot 1)\\
&=\p_\A^{(k)}L\ot r.
\end{align*}

\bigskip

{\em Step 5:\ In the notation of Step 3, $\chi=r_\B$.  That is,
$\Ctd_k(\B)=R_\B$.}

{\em Proof:} For any $a\in\B$, Step 4 shows that
\begin{align*}
\chi\left(a_{(1)}\p_\A^{(\ell)}L\ot 1\right)&=a_{(1)}\chi\left(\p_\A^{(\ell)}L\ot 1\right)\\
&=a_{(1)}\p_\A^{(\ell)}L\ot r\\
&=r_\B\left(a_{(1)}\p_\A^{(\ell)} L\ot 1\right).
\end{align*}
Since $\B$ is spanned by elements of the form
$a_{(1)}\p_\A^{(\ell)}L\ot 1$ (Step 2), we see that $\chi=r_\B$.

\bigskip

{\em Step 6:\ Two twisted loop algebras of $\A$ are $\R$-isomorphic
if and only if they are $k$-isomorphic.  That is, (i) and (ii) are
equivalent.}

{\em Proof:} Clearly
$R_\B\subseteq_\B\subseteq\Ctd_R(\B)\subseteq\Ctd_k(\B)$.  By Step
5, these inclusions are equalities.  Thus (i) and (ii) are
equivalent by Corollary \ref{3.27}.

\bigskip

From the eigenspace decomposition of $\widehat{\A}$, we see that the
left multiplication operator $L\ot 1_{(1)}:\
\widehat{\A}\rightarrow\widehat{\A}$ is injective, and it follows
easily that $Z(\B)=0$ for any twisted loop algebra $\B$ of $\A$.
Thus by Corollary \ref{uceiso}, (ii) is equivalent to (iv).

\bigskip

It remains to check that (i), (iii), and (iv) are equivalent. By
Theorem \ref{forms}, Proposition \ref{limit}, and Proposition
\ref{auto4},
 the $\R$-isomorphism classes of
${\wh\cS}/\R$-forms of the $\R$-conformal superalgebra
 $\A\ot\R$ are parametrized by \\$H^1_{\ct}\big(\wh{\Z},
\SL_2({\wh S} )\times \SL_2 (k)\big)/\pm(I,I)\big).$

 To simplify the notation, we write $G= \SL_2(\widehat{S}  ) \times \SL_2(k),$ $\overline{G}= \big(\SL_2(\widehat{S} )
 \times \SL_2(k)\big)/\pm(I,I)
 ,$ and consider the exact sequence of (continuous) $\widehat{\Z}$-groups
$$
1 \to \Z/2\Z \to G\to \overline{G} \to 1.
$$
where $\Z/2\Z$ is identified with the subgroup of $G$ generated by
$(-I,-I).$ This yields the exact sequence of pointed sets
\begin{equation}\label{coho4}
H_{\ct}^1(\wh{\Z},\Z/2\Z) \os\alpha  \to H_{ct}^1(\wh{\Z},G)
\os\beta \to H_{ct}^1(\wh{\Z},\overline{G}) \os{\p} \to
H_{\ct}^2(\wh{\Z},\Z/2\Z).
\end{equation}
whose individual terms can be understood as follows:

\begin{description}
  \item  {(a)}   $H^1_{\ct}(\wh{\Z},\Z/2\Z) = \,\Hom_{\ct}(\wh{\Z},\Z/2\Z)  \simeq \Z/2\Z. $
\newline A representative cocycle of the nonzero element of this
$H^1$ is the (unique and continuous) map $\wh{\Z} \to \Z/2\Z$ such
that $1 \mapsto (-I,-I).$

\item  {(b)} $H^2_{\ct}(\wh{\Z},\Z/2\Z) =0.$  This could be
checked by brute force in terms of cocycles.  An abstract argument
is as follows. Since $k$ is of characteristic $0$ we have $\Z/2\Z
\simeq \mu_2.$ Now $H^2_{\ct}(\wh{\Z},\mu  _2)$ is part (in fact all
of) the 2-torsion of the (algebraic) Brauer group
$H^2_{\et}(R,\bG_m).$ Because $R$ is of cohomological dimension $1,$
this Brauer group vanishes.

\item  {(c)} $H^1_{\ct}\big(\wh{\Z},\SL_2(\widehat{S}  )\big) =
1.$ Indeed this $H^1_{\ct}$ is the part of $H^1_{\et}(R,\SL_2)$
corresponding to the isomorphisms classes of $R$--torsors under
$\SL_2$ that become trivial over $\widehat{S}.$\footnote{In fact all
of $H^1_{\et}(R,\SL_2)$: Every $R$-torsor under $\SL_2$ is
isotrivial.} As observed in Remark \ref{comparison}
$H^1_{\et}(R,\SL_2)$ vanishes. \footnote{One can avoid the general
considerations of \cite{pianzola} in the present case by the
following direct argument.  The exact sequence of $R$-groups $ 1\to
\SL_2\to \bGL_2\os{\text{\rm det}}\to \,\bG_m \to 1 $ yields
$$
{\bf GL}_2(R) \os{\det}\to \, R^\times \to H^1_{\et}(R,\SL_2)\to
H^1_{\et}(R,\bGL_2).
$$
Since the map $\det$ is surjective, the map $H^1_{\et}(R,\SL_2)\to
H^1_{\et}(R,\bGL_2)$ has trivial kernel. On the other hand
$H^1_{\et}(R,\bGL_2)\simeq H^1_{\Zar}(R,\bGL_2) =1$ (the first
equality by Grothendieck-Hilbert 90, and the last since all rank 2
projective modules over $R$ are free; because $R$ is a principal
ideal domain).}
\end{description}

From (c), we immediately obtain

\begin{description}
\item{(d)} {\it The canonical map}
$H^1_{\ct}\big(\wh{\Z},\SL_2(k)\big) \to H^1_{\ct}(\wh{\Z},G)$ {\it
is bijective}.
\end{description}

Finally, since $\wh{\Z}$ acts trivially on $\SL_2(k)$ we have

\begin{description}
\item  {(e)} $H^1_{\ct}\big(\wh{\Z}, \SL_2(k)\big) \simeq$ \{conjugacy
classes of elements of finite order in $\SL_2(k)\}.$
\end{description}
By (b) and (e), we have a surjective map
\begin{align*}
\beta  :\ \ &\{ \;\hbox{\rm conjugacy classes of elements of finite
order
in} \; \SL_2(k)\}\\
 &\quad\quad\quad\quad\quad\quad\quad\quad\quad\quad\quad\quad\quad\quad\quad\quad\quad\quad\quad\quad\quad\quad \longrightarrow H^1_{\ct}(\wh{\Z},\overline{G}).
\end{align*}
Tracing through the various definitions, we see that the explicit
nature of the map $\beta$ is as follows: let $\theta \in \SL_2(k)$
be of finite order. Define a cocycle $u_\theta \in
Z^1(\wh{\Z},\overline{G})$ by $u_\theta (n) = \ol{(1,\theta ^n)}$
where $^- \, : \, G \to \overline{G}$ is the canonical map. Then
$\beta$ maps the conjugacy class of $\theta  $ to the class of
$u_\theta $ in $H^1_{\ct}(\wh{\Z},\overline{G}).$

It remains to show that $u_\theta  $ and $u_\sigma  $ are
cohomologous if and only if $\theta  $ is conjugate to $\pm\sigma .$
If $[u_\theta  ] = [u_\sigma  ]$ there exists $(x,\tau  )\in G =
\SL_2(\widehat{S}  )\times \SL_2(k)$ such that
$$
\ol{(x,\tau  )}\,^{-1} u_\theta  (n) \;^n\,\ol{(x,\tau  )} =
u_\sigma (n)
$$
for all $n\in \Z\subseteq \wh{\Z}.$  In particular, for $n = 1$ we
get
$$
\ol{({x^{-1}}\, {^{1}x},\tau  ^{-1}\theta  \tau  )} = \ol {(1,\sigma
)}
$$
which forces either
$$
\begin{array}{rlcrl}
x^{-1}\,^1x &= 1 &\q\text{\rm and} &\q \tau  ^{-1}\theta  \tau
&=\sigma\\
\text{\rm or}\q x^{-1}\,^1x &= -1 &\q\text{\rm and} &\q \tau
^{-1}\theta  \tau   &= -\sigma  .
\end{array}
$$
Thus $\theta$ is conjugate to either $\sigma$ or $- \sigma.$ The
converse is obvious given that the element
$\begin{pmatrix} t^{1/2} & 0 \\
0 & t^{1/2} \end{pmatrix} \in \SL_2(\widehat{S}^\times)$ satisfies
$x^{-1}\,^1x = -1.$

We have thus shown that $H^1_{\ct}(\wh{\Z},\overline{G})$ is in
bijective correspondence with the conjugacy classes of elements of
finite order in ${\bf PGL}_2(k).$ This completes the proof of the
theorem. \qed

\end{proof}

\bigskip

The grading operator $L$ is stable under all automorphisms of the
$N=2$ and $N=4$ conformal superalgebras, so $L\ot 1$ is an element
of every twisted loop algebra of these conformal superalgebras. By
considering the $n$-products of elements with $L\ot 1$, it is easy
to verify that every twisted loop algebra of the $N=2$ and $N=4$
conformal superalgebras is centreless.  They each admit a (unique)
one-dimensional universal central extension, as was previously shown
by one of the authors \cite{KvdL}.  Using Corollary \ref{uceiso}, we
see that Theorems \ref{4.20} and \ref{3.63} actually give a
parametrization of the $k$-isomorphism classes of universal central
extensions of twisted loop algebras of the $N=2$ and $N=4$ conformal
superalgebras.  Summarizing:

\begin{corollary}
There are exactly two $\C$-isomorphism classes of twisted loop
algebras based on the $N=2$ conformal superalgebra, and  infinitely
many $\C$-isomorphism classes of twisted loop algebras based on the
$N=4$ conformal superalgebra. The explicit automorphisms giving the
distinct isomorphism classes are the identity map and the
automorphism $\omega$ in the $N=2$ case; they are parametrized by
the conjugacy classes of elements of finite order in ${\bf
PGL}_2(\C)$ in the $N=4$ case.

Furthermore, two of these twisted loop algebras are isomorphic if and
only if their universal central extensions are isomorphic.\qed
\end{corollary}

\begin{remark} {\rm As explained in the Introduction, the superconformal
algebras  in Schwimmer and Seiberg's original work are obtained as
formal distribution algebras of the twisted loop algebras of the $N$
Lie conformal superalgebras. Since isomorphic twisted loop algebras
lead to isomorphic superconformal algebras, our work shows that in
the $N = 2, 4$ case  there can be no more superconformal algebras
than those listed by Schwimmer and Seiberg.}
\end{remark}

\begin{remark} {\rm Our methods work equally well for the other $N$-conformal superalgebras: The isomorphism classes of loop algebras based on an $N$-conformal superalgebra $\A$ are parametrized by
$$H^1\big(\widehat{\cS}  /\R,\,\autfun(\A)\big) \simeq
H_{\ct}^1\big(\pi_1(R),\,\autfun\,(\A)(\widehat{\cS}  )\big) \simeq
H_{\ct}^1\big(\widehat{\Z} ,\,\autfun\,(\A)(\widehat{\cS}  )\big)$$

For $N = 0$, the group $\autfun\,(\A)(\widehat{\cS}  )$ is trivial
and the only loop algebra is the affinization of $\A.$ For $N = 1$,
we have $\autfun\,(\A)(\widehat{\cS}  ) \simeq \Z/2\Z.$ There are
thus two non-isomorphic loop algebras (Ramond and Neveu-Schwarz).
For $N = 3$, the group $\autfun\,(\A)(\widehat{\cS}  )$ would appear
to be ${\rm \bf O}_3(\widehat{\cS}).$\footnote{The calculation is
delicate, just as in the $N =4$ case. We have not checked the
details thoroughly.} Under this assumption, the cohomology will
yield two non-isomorphic twisted loop algebras--again Ramond and
Neveu-Schwarz.}

\end{remark}

\end{document}